%
%
%
%


\magnification=1200
\pretolerance=500 \tolerance=1000 \brokenpenalty=5000
\hsize=12.5cm   
\vsize=19cm
\hoffset=0.4cm
\voffset=1cm   
\parskip3pt plus 1pt
\parindent=0.6cm
\let\sl=\it
\def\\{\hfil\break}
\def\emph#1{{\it #1}}
\def\frac#1#2{{#1\over #2}}


\font\seventeenbf=cmbx10 at 17.28pt

\font\twelvebf=cmbx10 at 12pt
\font\eightbf=cmbx8
\font\sixbf=cmbx6

\font\eighti=cmmi8
\font\sixi=cmmi6

\font\eightrm=cmr8
\font\sixrm=cmr6

\font\eightsy=cmsy8
\font\sixsy=cmsy6

\font\eightit=cmti8
\font\eighttt=cmtt8
\font\eightsl=cmsl8

\font\seventeenbsy=cmbsy10 at 17.28pt

\font\twelvebsy=cmbsy10 at 12pt
\font\tenbsy=cmbsy10
\font\eightbsy=cmbsy8
\font\sevenbsy=cmbsy7
\font\sixbsy=cmbsy6
\font\fivebsy=cmbsy5

\font\tenmsa=msam10

\font\sevenmsa=msam7
\font\fivemsa=msam5
\newfam\msafam
  \textfont\msafam=\tenmsa
  \scriptfont\msafam=\sevenmsa
  \scriptscriptfont\msafam=\fivemsa

\font\tenmsb=msbm10
\font\eightmsb=msbm8
\font\sevenmsb=msbm7
\font\fivemsb=msbm5
\newfam\msbfam
  \textfont\msbfam=\tenmsb
  \scriptfont\msbfam=\sevenmsb
  \scriptscriptfont\msbfam=\fivemsb
\def\Bbb{\fam\msbfam\tenmsb}

\font\tenCal=eusm10
\font\sevenCal=eusm7
\font\fiveCal=eusm5
\newfam\Calfam
  \textfont\Calfam=\tenCal
  \scriptfont\Calfam=\sevenCal
  \scriptscriptfont\Calfam=\fiveCal
\def\Cal{\fam\Calfam\tenCal}

\font\teneuf=eusm10
\font\teneuf=eufm10
\font\seveneuf=eufm7
\font\fiveeuf=eufm5
\newfam\euffam
  \textfont\euffam=\teneuf
  \scriptfont\euffam=\seveneuf
  \scriptscriptfont\euffam=\fiveeuf

\font\seventeenbfit=cmmib10 at 17.28pt

\font\twelvebfit=cmmib10 at 12pt
\font\tenbfit=cmmib10
\font\eightbfit=cmmib8
\font\sevenbfit=cmmib7
\font\sixbfit=cmmib6
\font\fivebfit=cmmib5
\newfam\bfitfam
  \textfont\bfitfam=\tenbfit
  \scriptfont\bfitfam=\sevenbfit
  \scriptscriptfont\bfitfam=\fivebfit


\catcode`\@=11
\def\eightpoint{%
  \textfont0=\eightrm \scriptfont0=\sixrm \scriptscriptfont0=\fiverm
  \def\rm{\fam\z@\eightrm}%
  \textfont1=\eighti \scriptfont1=\sixi \scriptscriptfont1=\fivei
  \def\oldstyle{\fam\@ne\eighti}%
  \textfont2=\eightsy \scriptfont2=\sixsy \scriptscriptfont2=\fivesy
  \textfont\itfam=\eightit
  \def\it{\fam\itfam\eightit}%
  \textfont\slfam=\eightsl
  \def\sl{\fam\slfam\eightsl}%
  \textfont\bffam=\eightbf \scriptfont\bffam=\sixbf
  \scriptscriptfont\bffam=\fivebf
  \def\bf{\fam\bffam\eightbf}%
  \textfont\ttfam=\eighttt
  \def\tt{\fam\ttfam\eighttt}%
  \textfont\msbfam=\eightmsb
  \def\Bbb{\fam\msbfam\eightmsb}%
  \abovedisplayskip=9pt plus 2pt minus 6pt
  \abovedisplayshortskip=0pt plus 2pt
  \belowdisplayskip=9pt plus 2pt minus 6pt
  \belowdisplayshortskip=5pt plus 2pt minus 3pt
  \smallskipamount=2pt plus 1pt minus 1pt
  \medskipamount=4pt plus 2pt minus 1pt
  \bigskipamount=9pt plus 3pt minus 3pt
  \normalbaselineskip=9pt
  \setbox\strutbox=\hbox{\vrule height7pt depth2pt width0pt}%
  \let\bigf@ntpc=\eightrm \let\smallf@ntpc=\sixrm
  \normalbaselines\rm}
\catcode`\@=12

\def\eightpointbf{%
 \textfont0=\eightbf   \scriptfont0=\sixbf   \scriptscriptfont0=\fivebf
 \textfont1=\eightbfit \scriptfont1=\sixbfit \scriptscriptfont1=\fivebfit
 \textfont2=\eightbsy  \scriptfont2=\sixbsy  \scriptscriptfont2=\fivebsy
 \eightbf
 \baselineskip=10pt}

\def\tenpointbf{%
 \textfont0=\tenbf   \scriptfont0=\sevenbf   \scriptscriptfont0=\fivebf
 \textfont1=\tenbfit \scriptfont1=\sevenbfit \scriptscriptfont1=\fivebfit
 \textfont2=\tenbsy  \scriptfont2=\sevenbsy  \scriptscriptfont2=\fivebsy
 \tenbf}
        
\def\twelvepointbf{%
 \textfont0=\twelvebf   \scriptfont0=\eightbf   \scriptscriptfont0=\sixbf
 \textfont1=\twelvebfit \scriptfont1=\eightbfit \scriptscriptfont1=\sixbfit
 \textfont2=\twelvebsy  \scriptfont2=\eightbsy  \scriptscriptfont2=\sixbsy
 \twelvebf
 \baselineskip=14.4pt}

\def\seventeenpointbf{%
 \textfont0=\seventeenbf  \scriptfont0=\twelvebf  \scriptscriptfont0=\eightbf
 \textfont1=\seventeenbfit\scriptfont1=\twelvebfit\scriptscriptfont1=\eightbfit
 \textfont2=\seventeenbsy \scriptfont2=\twelvebsy \scriptscriptfont2=\eightbsy
 \seventeenbf
 \baselineskip=20.736pt}
 

\newdimen\srdim \srdim=\hsize
\newdimen\irdim \irdim=\hsize
\def\NOSECTREF#1{\noindent\hbox to \srdim{\null\dotfill ???(#1)}}
\def\SECTREF#1{\noindent\hbox to \srdim{\csname REF\romannumeral#1\endcsname}}
\def\INDREF#1{\noindent\hbox to \irdim{\csname IND\romannumeral#1\endcsname}}
\newlinechar=`\^^J
\def\openauxfile{
  \immediate\openin1\jobname.
  \ifeof1
  \let\sectref=\NOSECTREF \let\indref=\NOSECTREF
  \else
  \input \jobname.aux
  \let\sectref=\SECTREF \let\indref=\INDREF
  \fi
  \immediate\openout1=\jobname.aux
  \let\END=\end \def\end{\immediate\closeout1\END}}
        
\newbox\titlebox   \setbox\titlebox\hbox{\hfil}
\newbox\sectionbox \setbox\sectionbox\hbox{\hfil}
\def\folio{\ifnum\pageno=1 \hfil \else \ifodd\pageno
           \hfil {\eightpoint\copy\sectionbox\kern8mm\number\pageno}\else
           {\eightpoint\number\pageno\kern8mm\copy\titlebox}\hfil \fi\fi}
\footline={\hfil}
\headline={\folio}

\def\titlerunning#1{\setbox\titlebox\hbox{\eightpoint #1}}
\def\title#1{\noindent\hfil$\smash{\hbox{\seventeenpointbf #1}}$\hfil
             \titlerunning{#1}\medskip}

\newcount\numbersection \numbersection=-1
\def\sectionrunning#1{\setbox\sectionbox\hbox{\eightpoint #1}
  \immediate\write1{\string\def \string\REF 
      \romannumeral\numbersection \string{%
      \noexpand#1 \string\dotfill \space \number\pageno \string}}}
\def\section#1{%
  \par\vskip0.666cm\penalty -100
  \vbox{\baselineskip=14.4pt\noindent{{\twelvepointbf #1}}}
  \vskip2pt
  \penalty 500
  \advance\numbersection by 1
  \sectionrunning{#1}}

\def\subsection#1{%
  \par\vskip0.5cm\penalty -100
  \vbox{\noindent{{\tenpointbf #1}}}
  \vskip1pt
  \penalty 500}

\newcount\numberindex \numberindex=0  
\def\index#1#2{%
  \advance\numberindex by 1
  \immediate\write1{\string\def \string\IND #1%
     \romannumeral\numberindex \string{%
     \noexpand#2 \string\dotfill \space \string\S \number\numbersection, 
     p.\string\ \space\number\pageno \string}}}

\newdimen\itemindent \itemindent=\parindent

\def\item#1{\par\noindent\hangindent\itemindent%
            \rlap{#1}\kern\itemindent\ignorespaces}
\def\itemitem#1{\par\noindent\hangindent2\itemindent%
            \kern\itemindent\rlap{#1}\kern\itemindent\ignorespaces}
\def\itemitemitem#1{\par\noindent\hangindent3\itemindent%
            \kern2\itemindent\rlap{#1}\kern\itemindent\ignorespaces}

\long\def\claim#1|#2\endclaim{\par\vskip 5pt\noindent 
{\tenpointbf #1.}\ {\sl #2}\par\vskip 5pt}

\def\proof{\noindent{\sl Proof}}

\def\today{\ifcase\month\or
January\or February\or March\or April\or May\or June\or July\or August\or
September\or October\or November\or December\fi \space\number\day,
\number\year}

\catcode`\@=11
\newcount\@tempcnta \newcount\@tempcntb 
\def\timeofday{{%
\@tempcnta=\time \divide\@tempcnta by 60 \@tempcntb=\@tempcnta
\multiply\@tempcntb by -60 \advance\@tempcntb by \time
\ifnum\@tempcntb > 9 \number\@tempcnta:\number\@tempcntb
  \else\number\@tempcnta:0\number\@tempcntb\fi}}
\catcode`\@=12

\def\bibitem#1&#2&#3&#4&%
{\hangindent=1.66cm\hangafter=1
\noindent\rlap{\hbox{\eightpointbf #1}}\kern1.66cm{\rm #2}{\sl #3}{\rm #4.}} 


\def\bC{{\Bbb C}}

\def\bQ{{\Bbb Q}}
\def\bR{{\Bbb R}}

\def\bZ{{\Bbb Z}}
\def\bone{{\mathchoice {\rm 1\mskip-4mu l} {\rm 1\mskip-4mu l}
{\rm 1\mskip-4.5mu l} {\rm 1\mskip-5mu l}}}


\def\cC{{\Cal C}}
\def\cE{{\Cal E}}

\def\cO{{\Cal O}}
\def\cH{{\Cal H}}


\def\square{{\hfill \hbox{
\vrule height 1.453ex  width 0.093ex  depth 0ex
\vrule height 1.5ex  width 1.3ex  depth -1.407ex\kern-0.1ex
\vrule height 1.453ex  width 0.093ex  depth 0ex\kern-1.35ex
\vrule height 0.093ex  width 1.3ex  depth 0ex}}}
\def\qed{\phantom{$\quad$}\hfill$\square$\medskip}
\def\hexnbr#1{\ifnum#1<10 \number#1\else
 \ifnum#1=10 A\else\ifnum#1=11 B\else\ifnum#1=12 C\else
 \ifnum#1=13 D\else\ifnum#1=14 E\else\ifnum#1=15 F\fi\fi\fi\fi\fi\fi\fi}
\def\msatype{\hexnbr\msafam}
\def\msbtype{\hexnbr\msbfam}
\mathchardef\restriction="3\msatype16   
\mathchardef\compact="3\msatype62
\mathchardef\smallsetminus="2\msbtype72   \let\ssm\smallsetminus
\mathchardef\subsetneq="3\msbtype28
\mathchardef\supsetneq="3\msbtype29
\mathchardef\leqslant="3\msatype36   \let\le\leqslant
\mathchardef\geqslant="3\msatype3E   \let\ge\geqslant
\mathchardef\ltimes="2\msbtype6E
\mathchardef\rtimes="2\msbtype6F

\let\ol=\overline

\let\wh=\widehat
\let\text=\hbox
\def\build#1|#2|#3|{\mathrel{\mathop{\null#1}\limits^{#2}_{#3}}}


\def\Re{\mathop{\rm Re}\nolimits}

\def\mod{\mathop{\rm mod}\nolimits} 

\def\Vol{\mathop{\rm Vol}\nolimits}
\def\MA{\mathop{\rm MA}\nolimits}
\def\PSH{\mathop{\rm PSH}\nolimits}

\def\exph{\mathop{\rm exph}\nolimits}
\def\logh{\mathop{\rm logh}\nolimits}

\def\dbar{{\overline\partial}}
\def\ddbar{{\partial\overline\partial}}

\def\alg{{\rm alg}}

\def\sing{{\rm sing}}
\def\loc{{\rm loc}}
\def\can{{\rm can}}

\long\def\InsertFig#1 #2 #3 #4\EndFig{\par
\hbox{\hskip #1mm$\vbox to#2mm{\vfil\special{" 
(/home/demailly/psinputs/grlib.ps) run
#3}}#4$}}
\long\def\LabelTeX#1 #2 #3\ELTX{\rlap{\kern#1mm\raise#2mm\hbox{#3}}}


\openauxfile

\title{Regularity of plurisubharmonic upper}
\smallskip
\title{envelopes in big cohomology classes}
\titlerunning{Regularity of plurisubharmonic upper envelopes in big 
cohomology classes}
\vskip2mm
{\leftskip=0.3cm
{\parindent=0cm
\hbox{\vbox{
{\bf Robert Berman}\\
Department of Mathematics,\\
Chalmers University of Technology,\\
Eklandag.\ 86,\\
SE-412 96 G\"oteborg, Sweden\\
{\it e-mail\/}: {\tt robertb@chalmers.se}}
\kern-6.4cm
\vbox{
{\bf Jean-Pierre Demailly}\\
Universit\'e de Grenoble I,\\
D\'epartement de Math\'ematiques,\\
Institut Fourier, BP 74,\\
38402 Saint-Martin d'H\`eres, France\\
{\tt demailly@fourier.ujf-grenoble.fr}
}}}\medskip}

\line{\hfill\it
dedicated to Professor Oleg Viro for his deep contributions to mathematics}
\bigskip

\noindent
{\bf Abstract.} The goal of this work is to prove the regularity of
certain quasi-plurisubharmonic upper envelopes. Such envelopes
appear in a natural way in the construction of hermitian metrics with 
minimal singularities on a big line bundle over a compact 
complex manifold. We prove that the complex Hessian forms of these 
envelopes are locally bounded outside an analytic set of singularities.
It~is furthermore shown that a parametrized version of this result
yields a priori inequalities for the solution of the Dirichlet problem
for a degenerate Monge-Amp\`ere operator$\,$; applications to geodesics
in the space of K\"ahler metrics are discussed.
A similar technique provides a logarithmic modulus of continuity for
Tsuji's ``supercanonical'' metrics, which generalize a well-known
construction of Narasimhan-Simha.
\smallskip

\noindent
{\bf R\'esum\'e.} Le but de ce travail est de d\'emontrer la
r\'egularit\'e de certaines enveloppes sup\'erieures de fonctions
quasi-plurisousharmoniques.  De telles enveloppes apparaissent
naturellement dans la construction des m\'etriques hermi\-tiennes \`a
singularit\'es minimales sur un fibr\'e en droites gros au dessus
d'une vari\'et\'e complexe compacte. Nous montrons que ces enveloppes
poss\`edent un Hessien complexe localement born\'e en dehors d'un
ensemble analytique de singularit\'es$\,$; par ailleurs, une version avec
param\`etres de ce r\'esultat permet d'obtenir des in\'egalit\'es a priori
pour la solution du probl\`eme de Dirichlet relatif \`a un op\'erateur de
Monge-Amp\`ere d\'eg\'en\'er\'e. Une technique similaire fournit
un module de continuit\'e logarithmique pour les m\'etriques
``super-canoniques'' de Tsuji, lesquelles g\'en\'eralisent une cons\-truction
bien connue de Narasimhan-Simha.
\smallskip

\noindent
{\bf Key words.} Plurisubharmonic function, upper envelope, 
hermitian line bundle, singular metric, logarithmic poles,
Legendre-Kiselman transform, pseudo-effective cone, volume,
Monge-Amp\`ere measure, supercanonical metric, Ohsawa-Takegoshi theorem.
\smallskip

\noindent
{\bf Mots-cl\'es.} Fonction plurisubharmonique, enveloppe sup\'erieure, 
fibr\'e en droites hermitien, m\'etrique singuli\`ere, p\^oles logarithmiques,
transform\'ee de Legendre-Kiselman, c\^one pseudo-effectif, volume,
mesure de Monge-Amp\`ere, m\'etrique super-canonique, 
th\'eor\`eme de Ohsawa-Takegoshi.
\smallskip

\noindent
{\bf AMS Classification.} 32F07, 32J25, 14B05, 14C30
\vfill\eject

\section{1. Main regularity theorem}

Let $X$ be a compact complex manifold and $\omega$ a hermitian metric
on~$X$, viewed as a smooth positive $(1,1)$-form. As usual we put
$d^c=\smash{1\over 4i\pi}(\partial-\dbar)$ so that 
$dd^c=\smash{1\over 2i\pi}\ddbar$. Consider 
the $dd^c$-cohomology class $\{\alpha\}$ of a smooth real \hbox{$d$-closed}
form $\alpha$ of type $(1,1)$ on $X$ [$\,$in general, one has to consider
the Bott-Chern cohomology group for which boundaries are $dd^c$-exact
$(1,1)$-forms $dd^c\varphi$, but in the case $X$ is K\"ahler, this group
is isomorphic to the Dolbeault cohomology group $H^{1,1}(X)\,$].
Recall that a function $\psi$ is said to be quasi-plurisubharmonic 
(or quasi-psh) if and only $idd^c\psi$ is locally bounded from below,
or equivalently, if it can be written locally as 
a sum $\psi=\varphi+u$ of a psh function $\varphi$ and a
smooth function~$u$. More precisely, it is said to be 
$\alpha$-plurisubharmonic (or $\alpha$-psh) if $\alpha+dd^c\psi\ge 0$.
We~denote by $\PSH(X,\alpha)$ the set of $\alpha$-psh functions on~$X$.

\claim (1.1) Definition|The class $\{\alpha\}\in H^{1,1}(X,\bR)$ is
said to be {\rm pseudo-effective} if it contains a closed 
$($semi-$)$positive current
$T=\alpha+dd^c\psi\ge 0$, and {\rm big} if it contains
a closed ``K\"ahler current'' $T=\alpha+dd^c\psi$ such that 
$T\ge\varepsilon\omega>0$ for some $\varepsilon>0$.
\endclaim

From now on in this section, we assume that $\{\alpha\}$ is \emph{big}.
We know by [Dem92] that we can then find $T_0\in\{\alpha\}$ of the form 
$$
T_0=\alpha+dd^c\psi_0\ge\varepsilon_0\omega\leqno(1.2)
$$ 
with a possibly slightly smaller $\varepsilon_0>0$ than the 
$\varepsilon$ in the definition, and $\psi_0$
a quasi-psh function with analytic singularities, i.e.\ locally
$$
\psi_0=c\log\sum|g_j|^2+u,~~~\hbox{where $c>0$, $u\in C^\infty$,
$g_j$ holomorphic}.\leqno(1.3)
$$
By [DP04], $X$ carries such a class $\{\alpha\}$ if and only if $X$
is in the Fujiki class $\cC$ of smooth varieties which are bimeromorphic
to compact K\"ahler manifolds. Our main result is

\claim (1.4) Theorem|Let $X$ be a compact complex manifold 
in the Fujiki class~$\cC$, and let
$\alpha$ be a smooth closed form of type $(1,1)$ on $X$ such that the
cohomology class $\{\alpha\}$ is big. Pick $T_0=\alpha+dd^c\psi_0
\in\{\alpha\}$ satisfying $(1.2)$ and $(1.3)$ for some hermitian metric
$\omega$ on~$X$, and let $Z_0$ be the analytic
set $Z_0=\psi_0^{-1}(-\infty)$. Then the upper envelope
$$
\varphi:=\sup\big\{\psi\le 0,~\psi~\hbox{$\alpha$-psh}\big\}
$$
is a quasi-plurisubharmonic function which has locally bounded second order 
derivatives $\partial^2\varphi/\partial z_j\partial\ol z_k$ on $X\ssm Z_0$,
and moreover, for suitable constants $C,\,B>0$, there is a global bound
$$
|dd^c\varphi|_\omega\le C(|\psi_0|+1)^2e^{B|\psi_0|}
$$
which explains how these derivatives blow up near $Z_0$. In particular
$\varphi$ is $C^{1,1-\delta}$ on $X\ssm Z_0$ for every $\delta>0$,
and the second derivatives $D^2\varphi$ are in $L^p_\loc(X\ssm Z_0)$
for every~$p>0$.
\endclaim

An important special case is the situation where we have a hermitian
line bundle $(L,h_L)$ and $\alpha=\Theta_{L,h_L}$, with the assumption that $L$
is big, i.e.\ that there exists a singular hermitian 
$h_0=h_L e^{-\psi_0}$ which has analytic singularities
and a curvature current 
$\Theta_{L,h_0}=\alpha+dd^c\psi_0\ge \varepsilon_0\omega$.
We then infer that the metric with minimal singularities
$h_{\min}=h_Le^{-\varphi}$ has the regularity properties prescribed
by Theorem~1.4 outside of the analytic set $Z_0=\psi_0^{-1}(-\infty)$.
In fact, [Ber07, Theorem~3.4~(a)] proves in this case the slightly
stronger result that $\varphi$ in $C^{1,1}$  on $X\ssm Z_0$ (using
the fact that $X$ is then Moishezon and that the total space of $L^*$
has a lot of holomorphic vector fields). The present approach is
by necessity different, since we can no longer rely on the existence
of vector fields when $X$ is not algebraic. Even then, our proof will 
be in fact somewhat simpler.

\proof. Notice that in order to get a quasi-psh function $\varphi$ we should
a priori replace $\varphi$ by its upper semi-continuous regularization
$\varphi^*(z)=\limsup_{\zeta \to z}\varphi(\zeta)$, but since 
$\varphi^*\le 0$ and $\varphi^*$ is $\alpha$-psh as well, $\psi=\varphi^*$
contributes to the envelope and therefore $\varphi=\varphi^*$.
Without loss of generality, after substracting a constant to
$\psi_0$, we may assume $\psi_0\le 0$. Then $\psi_0$ contributes to the
upper envelope and therefore $\varphi\ge\psi_0$. This already implies that
$\varphi$ is locally bounded on $X\ssm Z_0$. Following [Dem94], 
for every $\delta>0$, we consider the regularization operator
$$\psi\mapsto \rho_\delta\psi\leqno(1.5)$$
defined by $\rho_\delta\psi(z)=\Psi(z,\delta)$ and
$$\Psi(z,w)=\int_{\zeta\in T_{X,z}}\psi\big(\exph_z(w\zeta)\big)\,
\chi(|\zeta|^2)\,dV_\omega(\zeta),\qquad
(z,w)\in X\times\bC,\leqno(1.6)$$
where $\exph:T_X\to X$,~ $T_{X,z}\ni\zeta\mapsto\exph_z(\zeta)$
is the formal holomorphic part of the Taylor expansion of the exponential
map of the Chern connection on $T_X$ associated with the metric~$\omega$,
and $\chi:\bR\to\bR_+$ is a smooth function with support in $]-\infty,1]$ 
defined by 
$$
\chi(t)={C\over (1-t)^2}\exp{1\over t-1}\quad\hbox{for $t<1$},\qquad
\chi(t)=0\quad\hbox{for $t\ge 1$},
$$
with $C>0$ adjusted so that $\smash{\int_{|x|\le 1}}\chi(|x|^2)\,dx=1$ 
with respect to the Lebesgue measure $dx$ on $\bC^n$. Also, $dV_\omega(\zeta)$
denotes the standard hermitian Lebesgue measure on~$(T_X,\omega)$.
Clearly $\Psi(z,w)$ depends only on $|w|$. With the relevant change of 
notation, the estimates proved in sections 
3 and 4 of [Dem94] (see especially Theorem 4.1 and estimates (4.3), (4.5)
therein) show that if one assumes $\alpha+dd^c\psi\ge 0$,
there are constants $\delta_0,K>0$ such that for $(z,w)\in X\times\bC$
$$
\leqalignno{
&~~~~[0,\delta_0]\ni t\mapsto \Psi(z,t)+Kt^2\quad\hbox{is increasing,}&(1.7)\cr
\noalign{\vskip5pt}
&~~~~\alpha(z)+dd^c\Psi(z,w)\ge -A\lambda(z,|w|)|dz|^2
-K\big(|w|^2|dz|^2+|dz||dw|+|dw|^2),&(1.8)}
$$
where $A=\sup_{|\zeta|\le 1,|\xi|\le 1}\{-c_{jk\ell m}\zeta_j\ol\zeta_k
\xi_\ell\ol\xi_m\}$ is a bound for the negative part of the 
curvature tensor $(c_{jk\ell m})$ of $(T_X,\omega)$ and
$$
\lambda(z,t)={d\over d\log t}(\Psi(z,t)+Kt^2)
\build\longrightarrow~~||t\to 0_+|\nu(\psi,z)\qquad
\hbox{(Lelong number).}\leqno(1.9)
$$
In fact, this is clear from [Dem94] if $\alpha=0$, and otherwise
we simply apply the
above estimates (1.7--1.9) locally to $u+\psi$ where $u$ is a local
potential of $\alpha$ and then subtract the resulting regularization
$U(z,w)$ of $u$ which is such that 
$$
dd^c(U(z,w)-u(z))=O(|w|^2|dz|^2+|w||dz||dw|+|dw|^2)\leqno(1.10)
$$
because the
left hand side is smooth and $U(z,w)-u(z)=O(|w|^2)$. As a consequence,
the regularization operator $\rho_\delta$ transforms quasi-psh functions
into quasi-psh functions, while providing very good control on the complex
Hessian. We exploit this, again quite similarly as in [Dem94], by introducing 
the Kiselman-Legendre transform (cf.\ [Kis78, Kis94])
$$
\psi_{c,\delta}(z)=\inf_{t\in{}]0,\delta]}
\rho_t\psi(z)+Kt^2-K\delta^2-c\log{t\over\delta},\qquad
c>0,~\delta\in{}]0,\delta_0].\leqno(1.11)
$$
We need the following basic lower bound on the Hessian form.

\claim (1.12) Lemma|For all $c>0$ and $\delta\in{}]0,\delta_0]$ we have
$$
\alpha+dd^c\psi_{c,\delta}\ge -\big(A\min\big(c,\lambda(z,\delta)\big)
+K\delta^2\big)\omega.
$$
\endclaim

\proof\ {\it of lemma}. In general an
infimum $\inf_{\eta\in E}u(z,\eta)$ of psh functions 
$z\mapsto u(z,\eta)$ is not 
psh, but this is the case if $u(z,\eta)$ is psh with respect to
$(z,\eta)$ and $u(z,\eta)$ depends only on $\Re\eta$ -- in which case it
is actually a convex function of $\Re\eta$ -- this fundamental fact is
known as Kiselman's infimum principle. We apply it here by putting
$w=e^\eta$ and $t=|w|=e^{\Re\eta}$. At all points of 
$E_c(\psi)=\{z\in X\,;\;\nu(\psi,z)\ge c\}$ the infinimum occurring
in (1.11) is attained
at $t=0$. However, for $z\in X\ssm E_c(\psi)$ it is attained 
for $t=t_{\min}$ where 
$$\cases{
t_{\min}=\delta&if $\lambda(z,\delta)\le c$,\cr
t_{\min}<\delta&such that $c=\lambda(z,t_{\min})={d\over dt}(
\Psi(z,t)+Kt^2)_{t=t_{\min}}$ if $\lambda(z,\delta)>c$.\cr}
$$
In a neighborhood of such a point $z\in X\ssm E_c(\psi)$, 
the infimum coincides with the infimum taken for $t$
close to $t_{\min}$, and all functions involved have (modulo
addition of~$\alpha$) a Hessian form bounded below by 
$-(A\lambda(z,t_{\min})+K\delta^2)\omega$ by (1.8).
Since $\lambda(z,t_{\min})\le \min(c,\lambda(z,\delta))$, we get the
desired estimate on the dense open set $X\ssm E_c(\psi)$ by Kiselman's
infimum principle. However $\psi_{c,\delta}$ is quasi-psh on~$X$ and 
$E_c(\psi)$ is of measure zero, so the estimate is in fact valid on 
all of~$X$, in the sense of currents.\qed
\medskip

We now proceed to complete the proof or Theorem~1.4.
Lemma~1.12 implies the more brutal estimate
$$
\alpha+dd^c\psi_{c,\delta}\ge -(Ac+K\delta^2)\,\omega\qquad
\hbox{for $\delta\in{}]0,\delta_0]$}.\leqno(1.13)
$$
Consider the convex linear combination
$$
\theta={Ac+K\delta^2\over\varepsilon_0}\psi_0+\bigg
(1-{Ac+K\delta^2\over\varepsilon_0}\bigg)\varphi_{c,\delta}
$$
where $\varphi$ is the upper envelope of all $\alpha$-psh functions 
$\psi\le 0$. Since 
$\alpha+dd^c\varphi\ge 0$, (1.2) and (1.13) imply
$$
\alpha+dd^c\theta\ge (Ac+K\delta^2)\,\omega-
\bigg(1-{Ac+K\delta^2\over\varepsilon_0}\bigg)(Ac+K\delta^2)\,\omega\ge 0.
$$
Also $\varphi\le 0$ and therefore 
$\varphi_{c,\delta}\le\rho_\delta\varphi\le 0$ and $\theta\le 0$ likewise.
In particular $\theta$ contributes to the envelope and as a consequence we get
$\varphi\ge\theta$. Coming back to the definition of $\varphi_{c,\delta}$, we
infer that for every point $z\in X\ssm Z_0$ and every $\delta>0$, there exists 
$t\in{}]0,\delta]$ such that
$$
\eqalign{
\varphi(z)&\ge {Ac+K\delta^2\over\varepsilon_0}\psi_0(z)+
\bigg(1-{Ac+K\delta^2\over\varepsilon_0}\bigg)
(\rho_t\varphi(z)+Kt^2-K\delta^2-c\log t/\delta)\cr
&\ge {Ac+K\delta^2\over\varepsilon_0}\psi_0(z)+
(\rho_t\varphi(z)+Kt^2-K\delta^2-c\log t/\delta)\cr}
$$
(using the fact that the infimum is${}\le 0$ and reached for some
$t\in{}]0,\delta]$, as 
$t\mapsto \rho_t\varphi(z)$ is bounded for $z\in X\ssm Z_0$). Therefore we get
$$
\rho_t\varphi(z)+Kt^2\le 
\varphi(z)+K\delta^2- (Ac+K\delta^2)\varepsilon_0^{-1}\psi_0(z)
+c\log{t\over \delta}\;.\leqno(1.14)
$$
Since $t\mapsto \rho_t\varphi(z)+Kt^2$ is increasing and equal to $\varphi(z)$
for $t=0$, we infer that
$$
K\delta^2-(Ac+K\delta^2)\varepsilon_0^{-1}\psi_0(z)+c\log{t\over \delta}\ge 0,
$$
or equivalently, since $\psi_0\le 0$,
$$
t \ge \delta\,\exp\big(-(A+K\delta^2/c)\varepsilon_0^{-1}|\psi_0(z)|
-K\delta^2/c\big).
$$
Now, (1.14) implies the weaker estimate
$$
\rho_t\varphi(z)\le \varphi(z)+K\delta^2+(Ac+K\delta^2)
\varepsilon_0^{-1}|\psi_0(z)|,
$$
hence, by combining the last two inequalities, we get
$$
\eqalign{
&{\rho_t\varphi(z)-\varphi(z)\over t^2}\cr
&\qquad{}\le 
K\bigg(1+\Big({\textstyle Ac\over K\delta^2}+1\Big)
\varepsilon_0^{-1}|\psi_0(z)|\bigg)
\exp\bigg(2\Big(A+K{\textstyle \delta^2\over c}\Big)
\varepsilon_0^{-1}|\psi_0(z)|+2K{\textstyle \delta^2\over c}\bigg).\cr}
$$
We exploit this by letting $0<t\le\delta$ and $c$ tend to $0$, in such a way
that $Ac/K\delta^2$ converges to a positive limit~$\ell$ (if $A=0$,
just enlarge $A$ slightly and then let~$A\to 0$). In this way we 
get for every $\ell>0$
$$
\eqalign{
&\liminf_{t\to 0_+}~{\rho_t\varphi(z)-\varphi(z)\over t^2}\cr
&\qquad{}
\le K\big(1+(\ell+1)\varepsilon_0^{-1}|\psi_0(z)|\big)
\exp\Big(2A\big((1+\ell^{-1})\varepsilon_0^{-1}|\psi_0(z)|+
\ell^{-1}\big)\Big).\cr}
$$
The special (essentially optimal) choice 
$\ell=\varepsilon_0^{-1}|\psi_0(z)|+1$ yields
$$
\liminf_{t\to 0_+}~{\rho_t\varphi(z)-\varphi(z)\over t^2}\le
K(\varepsilon_0^{-1}|\psi_0(z)|+1)^2\exp
\big(2A(\varepsilon_0^{-1}|\psi_0(z)|+1)\big).\leqno(1.15)
$$
Now, putting as usual $\nu(\varphi,z,r)={1\over \pi^{n-1}r^{2n-2}/(n-1)!}
\int_{B(z,r)}\Delta\varphi(\zeta)\,d\zeta$, we infer from estimate~(4.5) 
of [Dem94] the Lelong-Jensen like inequality
$$
\leqalignno{
\rho_t\varphi(z)-\varphi(z)&=
\int_0^t{d\over d\tau}\Phi(z,\tau)\,d\tau\cr
&\ge\int_0^t{d\tau\over\tau}
\bigg(
\int_{B(0,1)}\nu(\varphi,z,\tau|\zeta|)\,\chi(|\zeta|^2)\,d\zeta
-O(\tau^2)\bigg)\cr
\noalign{\vskip5pt}
&\ge c(a)\,\nu(\varphi,z,at)-C_2t^2\qquad
\hbox{where $a<1$, $c(a)>0$ and $C_2\gg 1$}, \cr
\noalign{\vskip5pt}
&={c'(a)\over t^{2n-2}}
\int_{B(z,at)}\Delta\varphi(\zeta)\,d\zeta-C_2t^2&(1.16)\cr}
$$
where the third line is obtained
by integrating for $\tau\in[a^{1/2}t,t]$ and for $\zeta$ in the corona
$a^{1/2}<|\zeta|<a^{1/4}$
(here we assume that $\chi$ is taken to be decreasing with $\chi(t)>0$ 
for all $t<1$, and we compute the laplacian $\Delta$ in normalized
coordinates at $z$ given by $\zeta\mapsto\exph_z(\zeta)$). Hence
by Lebesgue's theorem on the existence almost everywhere of the
density of a positive measure (see e.g.\ [Rud66], 7.14), we find
$$
\lim_{t\to 0_+}{1\over t^2}\big(\rho_t\varphi(z)-\varphi(z)\big)\ge
c''(\Delta_\omega\varphi)_{\rm ac}(z)-C_2\quad \hbox{a.e.\ on $X$}\leqno(1.17)
$$
where the ac  subscript means the absolutely continuous part of the measure
$\Delta_\omega\varphi$. By combining (1.15) and (1.17) and using the
quasi-plurisubharmonicity of $\varphi$ we conclude that
$$
|dd^c\varphi|_\omega\le \Delta_\omega\varphi+C_3\le 
C\,(|\psi_0|+1)^2\,e^{2A\varepsilon_0^{-1}\psi_0(z)}\quad\hbox{a.e.\
on $X\ssm Z_0$}
$$
for some constant $C>0$. There cannot be any singular measure part $\mu$
in $\Delta_\omega\varphi$ either, since we now that the Lebesgue density 
would then be equal to $+\infty$ \hbox{$\mu$-a.e.}
([Rud66], 7.15), in contradiction with (1.15). This gives the
required estimates for the
complex derivatives $\partial^2\varphi/\partial z_j\partial\ol z_k$.
The other real derivatives $\partial^2\varphi/\partial x_i\partial x_j$
are obtained from $\Delta\varphi=\sum_k 
\partial^2\varphi/\partial z_k\partial\ol z_k$ via singular integral 
operators, and it is well-known that these operate boundedly
on $L^p$ for all $p<\infty$. Theorem~(1.4) follows.\qed

\claim (1.18) Remark|{\rm The proof gave us in fact the very explicit value
$B=2A\varepsilon_0^{-1}$, where $A$ is an upper bound of the negative
part of the curvature of $(T_X,\omega)$. The slightly more refined
estimates obtained in [Dem94] show that we could even replace 
$B$ by the possibly smaller constant 
$B_\eta=2(A'+\eta)\varepsilon_0^{-1}$ where
$$
A'=\sup_{|\zeta|=1,\,|\xi|=1,\,\zeta\perp\xi}
-c_{jk\ell m}\zeta_j\ol\zeta_k\xi_\ell\ol\xi_m,
$$
and the dependence of the other constants on~$\eta$ could then be
made explicit.}
\endclaim

\claim (1.19) Remark|{\rm In Theorem (1.4), one can replace the assumption
that $\alpha$ is smooth by the assumption that $\alpha$ has $L^\infty$
coefficients. In fact, we used the smoothness of $\alpha$ only as
a cheap argument to get the validity of estimate (1.10) for the local 
potentials $u$ of $\alpha$. However, the results of [Dem94] easily imply the
same estimates when $\alpha$ is $L^\infty$, as both $u$ and $-u$
are then quasi-psh; this follows e.g.\ from (1.8) applied with respect
to a smooth $\alpha_\infty$ and $\psi=\pm u$ 
if we observe that $\lambda(z,|w|)=O(|w|^2)$ when $|dd^c\psi|_\omega$ 
is bounded. Therefore, only the constant $K$ will be affected in the
proof.}
\endclaim

\section{2. Applications to volume and Monge-Amp\`ere measures }

Recall that the \emph{volume} of a big class $\{\alpha\}$ is defined,
in the work [Bou02] of S.~Boucksom, as
$$
\Vol(\{\alpha\})=\sup_{T}\int_{X\ssm\sing(T)}T^{n},\leqno(2.1)
$$
with $T$ ranging over all positive currents in the class $\{\alpha\}$
with \emph{analytic singularities}, whose locus is denoted by $\sing(T)$.
If the class is not big then the volume is defined to be zero. With
this definition, it is clear that $\{\alpha\}$ is big precisely
when $\Vol(\{\alpha\})>0$. 

Now fix a \emph{smooth} representative $\alpha$ in a \emph{pseudo-effective}
class $\{\alpha\}$. We then obtain a uniquely defined 
$\alpha$-plurisubharmonic function $\varphi=\psi_{\min}\ge 0$ with minimal 
singularities defined as in Theorem~(1.4) by
$$
\varphi:=\sup\big\{\psi\le 0,~\psi~\hbox{$\alpha$-psh}\big\}\,;
\leqno(2.2)
$$
notice that the supremum is non empty by our assumption that $\{\alpha\}$ is 
pseudo-effective. If $\{\alpha\}$ is big
and $\psi$ is
$\alpha$-psh and locally bounded in the complement of an analytic $Z\subset X$,
one can define the Monge-Amp\`ere measure $\MA_{\alpha}(\psi)$ by
$$
\MA_{\alpha}(\psi):=\bone_{X\ssm Z}(\alpha+dd^c\psi)^{n}
\leqno(2.3)
$$
as follows from the work of Bedford-Taylor [BT76, BT82].
In particular, if $\{\alpha\}$ is big, there is a well-defined 
positive measure on $\MA_{\alpha}(\varphi)=\MA_{\alpha}(\psi_{\min})$ on~$X\,$;
its total mass coincides with $\Vol(\{\alpha\})$, i.e.\
$$
\Vol(\{\alpha\})=\int_{X}\MA_{\alpha}(\varphi)
$$
(this follows from the comparison theorem and the fact that Monge-Amp\`ere
measures of locally bounded psh functions do not carry mass on analytic 
sets$\,$; see e.g.\ [BEGZ08]). Next, notice that in general
the $\alpha$-psh envelope $\varphi=\psi_{\min}$ corresponds \emph{canonically}
to $\alpha$, so we may associate to $\alpha$ the following subset of $X\,$:
$$
D=\{\varphi=0\}.
\leqno(2.4)
$$
Since $\varphi$ is upper semi-continuous, the set $D$ is compact. Moreover,
a simple application of the maximum principle shows that $\alpha\geq 0$
pointwise on $D$ (precisely as in Proposition 3.1 of [Ber07]: at any point 
$z_0$ where $\alpha$ is not semi-positive, we can find complex coordinates and
a small $\varepsilon>0$ such that $\varphi(z)-\varepsilon|z-z_0|^2$ 
is subharmonic near~$z_0$, hence $\varphi(z_0)<0$).
In particular, $\bone_{D}\alpha$ is a positive $(1,1)$-form on~$X$. From
Theorem~(1.4) we infer

\claim (2.5) Corollary|
Assume that $X$ is a K\"ahler manifold. For any smooth closed
form $\alpha$ of type $(1,1)$ in a pseudo-effective class and $\varphi\le 0$
the $\alpha$-psh upper envelope we have
$$
\MA_{\alpha}(\varphi)=\bone_{D}\alpha^n,\qquad
D=\{\varphi=0\},
\leqno(2.6)
$$
as measures on $X$ $($provided the left hand side is interpreted
as a suitable weak limit$)$ and
$$
\Vol(\{\alpha\})=\int_{D}\alpha^{n}\ge 0.
\leqno(2.7)
$$
 In particular, $\{\alpha\}$ is big if and only if $\int_{D}\alpha^{n}>0$.
\endclaim

\proof. Let $\omega$ be a K\"ahler metric on~$X$.
First assume that the class $\{\alpha\}$ is big and let $Z_0$ be
the singularity set of some strictly positive representative 
$\alpha+dd^c\psi_0\ge\varepsilon\omega$ with 
analytic singularities. By Theorem~(1.4),
$\alpha+dd^c\varphi$ is in $L_{\loc}^{\infty}(X\ssm Z_0)$. In
particular (see [Dem89]) 
the Monge-Amp\`ere measure $(\alpha+dd^c\varphi)^{n}$
has a locally bounded density on $X\ssm Z_0$ with respect 
to $\omega^{n}$. Hence,
it is enough to prove the identity (2.6) pointwise almost
enerywhere on $X$. To this end, one argues essentially as in [Ber07]
(where the class was assumed to be integral). First a well-known local 
argument based on the solution of the Dirichlet problem for $(dd^c)^n$
(see e.g.\ [BT76, BT82], and also Proposition 1.10 in [BB08]) proves that 
the Monge-Amp\`ere
measure $(\alpha+dd^c\varphi)^{n}$ of the envelope $\varphi$ vanishes
on the open set $(X\ssm Z_0)\ssm D$ (this only uses the fact that 
$\alpha$ has continuous potentials and the continuity of $\varphi$ on
$X\ssm Z_0$). Moreover, Theorem~(1.4) implies that $\varphi\in C^1(X\ssm Z_0)$
and
$$
\frac{\partial^2\varphi}{\partial x_{i}\partial x_{j}}\in L_{\loc}^{p}
\leqno(2.8)
$$
for any $p\in{}]1,\infty[$ and $i,j\in[1,2n]$. Even if this is slightly
weaker than the situation in [Ber07], where it was shown that one
can take $p=\infty$, the argument given in [Ber07] still goes
through. Indeed, by well-known properties of measurable sets, $D$ has
Lebesgue density $\lim_{r\to 0}\lambda(D\cap B(x,r))/\lambda(B(x,r)=1$
at almost every point $x\in D$, and since $\varphi=0$ on $D$, we conclude
that $\partial\varphi/\partial x_i=0$ at those points (if the
density is $1$, no open cone of vertex $x$ can be omitted and thus we can
approach $x$ from any direction by a sequence $x_\nu\to x$). 
But the first derivative is H\"older continuous on $D\ssm Z_0$, hence
$\partial\varphi/\partial x_i=0$ everywhere on $D\ssm Z_0$.
By repeating the argument for
$\partial\varphi/\partial x_i$ which has a derivative in $L^p$ ($L^1$ would
even be enough), we conclude from Lebesgue's theorem that
$\partial^2\varphi/\partial x_{i}\partial x_{j}=0$
a.e.\ on~$D\ssm Z_0$, hence $\alpha+dd^c\varphi=\alpha$ on~$D\ssm E$
where the set $E$ has measure zero with respect to $\omega^{n}$.
This proves formula (2.6) in the case of a big class. 

Finally, assume that $\{\alpha\}$ is pseudo-effective but not big.
For any given positive number $\varepsilon$ we 
let $\alpha_{\varepsilon}=\alpha+\varepsilon\omega$
and denote by $D_{\varepsilon}$ the corresponding set (2.4).
Clearly $\alpha_{\varepsilon}$ represents a big class. Moreover, by
the continuity of the volume function up to the boundary of the big
cone [Bou02]
$$
\Vol(\{\alpha_{\varepsilon}\})\rightarrow\Vol(\{\alpha\})\quad ({}=0)
\leqno(2.9)
$$
as $\varepsilon$ tends to zero. Now observe that 
$D\subset D_{\varepsilon}$
(there are more $(\alpha+\varepsilon\omega)$-psh 
functions than $\alpha$-psh functions and so 
$\varphi\le\varphi_\varepsilon\le 0\,$; clearly 
$\varphi_\varepsilon$ increases with $\varepsilon$ and 
$\varphi=\lim_{\varepsilon\to 0}\varphi_\varepsilon\,$; 
compare with Proposition 3.3 in [Ber07]). Therefore
$$
\int_{D}\alpha^{n}\leq\int_{D_{\varepsilon}}\alpha^{n}\leq\int_{D_{\varepsilon}}\alpha_{\varepsilon}^{n},
$$
where we used that $\alpha\leq\alpha_{\varepsilon}$
in the second step. Finally, since by the big case treated above,
the right hand side above is precisely $\Vol(\{\alpha_{\varepsilon}\})$,
letting $\varepsilon$ tend to zero and using (2.9) proves
that $\int_{D}\alpha^{n}=0=\Vol(\{\alpha\})$ [and that
$\MA_\alpha(\varphi)=0$ if we interpret it as the limit of
$\MA_{\alpha_\varepsilon}(\varphi_\varepsilon)$]. This concludes the proof.\qed

In the case when $\{\alpha\}$ is an integer class, i.e. when it is
the first Chern class $c_1(L)$ of a holomorphic line bundle $L$
over $X$, the result of the corollary was obtained in [Ber07]
under the additional assumption that $X$ be a \emph{projective} manifold --
it was conjectured there that the result was also valid for integral
classes over non-projective K\"ahler manifold.

\claim (2.10) Remark|{\rm
In particular, the corollary shows that, if $\{\alpha\}$ is big,
there is always an $\alpha$-plurisubharmonic function $\varphi$ with
minimal singularities such that $\MA_{\alpha}(\varphi)$ has a
$L^{\infty}$-density
with respect to $\omega^{n}$. This is a very useful fact when dealing with
big classes which are not K\"ahler (see for example [BBGZ09]).}
\endclaim

\section{3. Application to regularity of a boundary value problem and 
a variational principle}

In this section we will see how the main theorem may be intepreted
as a regularity result for (1) a free boundary value problem for the
Monge-Amp\`ere operator and (2) a variational principle. For simplicity
we only consider the case of a K\"ahler class.

\subsection{3.1. A free boundary value problem for the
Monge-Amp\`ere operator}

Let $(X,\omega)$ be a K\"ahler manifold. Given a function $f\in\cC^2(X)$
consider the following \emph{free} boundary value problem
$$
\left\{
\eqalign{
\MA_{\omega}(u)&= 0\qquad\hbox{on $\Omega$},\cr
u &= f \qquad\hbox{on $\partial\Omega$},\cr
du & = df\cr}\right.
$$
for a pair $(u,\Omega)$, where $u$ is an $\omega$-psh function
on $\overline{\Omega}$ which is in $\cC^1(\overline{\Omega})$
and $\Omega$ is an open set in $X$. We have used the notation $\partial\Omega:=\overline{\Omega}\ssm\Omega$,
but no regularity of the boundary is assumed. The reason why the
set $\Omega$ is assumed to be part of the solution is that, for 
a fixed $\Omega$,
the equations are overdetermined. Setting $u:=\varphi+f$ and $\Omega:=X\ssm D$
where $\varphi$ is the upper envelope with respect to $\alpha:=dd^cf+\omega$,
yields a solution. In fact, by Theorem~(1.4)
$u\in\cC^{1,1-\delta}(\overline{\Omega})$
for any $\delta>0$.

\subsection{3.2. A variational principle}

Fix a form $\alpha$ in a K\"ahler class $\{\alpha\}$, possessing continuous
potentials. Consider the following energy functional defined on the
convex space $\PSH(X,\alpha)\cap L^{\infty}$ of all $\alpha$-psh
functions which are bounded on $X$~:
$$
\cE[\psi]:=\frac{1}{n+1}\sum_{j=0}^{n}\int_{X}\psi(\alpha+dd^c\psi)^{j}\wedge\alpha^{n-j}
\leqno(3.2.1)
$$
This functional seems to first have appeared, independently, in the
work of Aubin and Mabuchi in K\"ahler-Einstein geometry (in the case
when $\alpha$ is a K\"ahler form). More geometrically, up to an additive
constant, $\cE$ can be defined as a primitive of the one
form on $\PSH(X,\alpha)\cap L^{\infty}$ defined by the measure valued
operator $\psi\mapsto \MA_{\alpha}(\psi)$. 
$$
$$
As shown in [BB08] (version 1) the following variational characterization
of the envelope $\varphi$ holds:

\claim (3.2.2) Proposition|
The functional 
$$
\psi\mapsto\cE[\psi]-\int_{X}\psi(\alpha+dd^c\psi)^{n}
$$
achieves its minimum value on the space $\PSH(X,\alpha)\cap L^{\infty}$
precisely when $\psi$ is equal to the envelope $\varphi$ 
$($defined with respect to $\alpha)$. Moreover, the minimum
is achieved only at $\varphi$, up to an additive constant.
\endclaim

Hence, the main theorem above can be interpreted as a regularity result
for the functions in $\PSH(X,\alpha)\cap L^{\infty}$ minimizing the
functional (3.2.1) in the case when $\alpha$ is assumed to have
$L^\infty_\loc$ coefficients.
More generally, a similar variational characterization of $\varphi$
can be given the case of a \emph{big} class $[\alpha]$ [BBGZ09].

\section{4. Degenerate Monge-Amp\`ere equations and geodesics 
in the space of K\"ahler metrics}

Assume that $(X,\omega)$ is a compact K\"ahler manifold and that $\Sigma$
is a Stein manifold with strictly pseudoconvex boundary, i.e. $\Sigma$ 
admits a smooth
strictly psh non-positive function $\eta_\Sigma$ which vanishes precisely
on~$\partial\Sigma$. The corresponding product manifold will be denoted
by $M:=\Sigma\times X$. By taking pull-backs, we identify $\eta_\Sigma$
with a function on $M$ and $\omega$ with a semi-positive form
on $M$. In this way, we obtain a K\"ahler form $\omega_M:=\omega+dd^c\eta_\Sigma$
on~$M$. Given a function $f$ on $M$ and a point $s$ in $\Sigma$ we
use the notation $f_s:=f(s,\cdot)$ for the induced function
on $X$.

Further, given a closed $(1,1)$ form $\alpha$ on $M$ with bounded coefficients
and a continuous function $f$ on $\partial M$, we define the upper envelope:
$$
\varphi_{\alpha,f}:=\sup\left\{ \psi:\,\psi\in\PSH(M,\alpha)\cap C^0(M),\,\,\psi_{\partial M}\leq f\right\}.\leqno(4.1)
$$
Note that when $\Sigma$ is a point and $f=0$ this definition coincides with
the one introduced in section~1. Also, when $F$ is a smooth function
on the whole of $M$, the obvious translation $\psi\mapsto\psi'=\psi-F$ yields 
the relation
$$
\varphi_{\beta,f-F}=\varphi_{\alpha,f}-F\quad\hbox{where}~~
\beta=\alpha+dd^cF.\leqno(4.2)
$$
The proof of the following lemma is a straightforward adaptation
of the proof of Bedford-Taylor [BT76] in the case when $M$ is a strictly
pseudoconvex domain in~$\bC^n$.

\claim (4.3) Lemma|
Let $\alpha$ be a closed real $(1,1)$-form on $M$ with bounded 
coefficients, such that $\alpha_{|\{s\}\times X}\ge\varepsilon_0\omega$ is 
positive definite for all $s\in\Sigma$.
Then the corresponding envelope $\varphi=\varphi_{\alpha,0}$ vanishes on 
the boundary of $M$ and is continuous on $M$. Moreover, 
$MA_{\alpha}(\varphi)$ vanishes in the interior of~$M$.
\endclaim

\proof. By (4.2) we have $\varphi_{\alpha,0}=\varphi_{\beta,0}+C
\eta_\Sigma$ where $\beta=\alpha+C dd^c\eta_\Sigma$ can be taken to be positive
definite on~$M$ for $C\gg 1$, as is easily seen from the Cauchy-Schwarz 
inequality and the
hypotheses on~$\alpha$. Therefore, we can assume without loss of generality that
$\alpha$ is positive definite on~$M$.
Since $0$ is a candidate for the supremum defining $\varphi$ it follows immediately
that $0\leq\varphi$ and hence $\varphi_{\partial M}=0$. To see that $\varphi$
is continuous on $\partial M$ (from the inside) take an arbitrary
candidate $\psi$ for the sup and observe that
$$
\psi\leq-C\eta_\Sigma
$$
for $C\gg 1$, independent of $\psi$. Indeed, since $dd^c\psi\geq-\alpha$
there is a large positive constant $C$ such that the function $\psi+C\eta_\Sigma$
is strictly plurisubharmonic on $\Sigma\times\{ x\}$ for all $x.$
Thus the inequality above follows from the maximum principle applied
to all slices $\Sigma\times\{ x\}$. All in all, taking the sup over
all such $\psi$ gives
$$
0\leq\varphi\leq-C\eta_\Sigma.
$$
But since $\eta_{\Sigma|\partial M}=0$ and $\eta_\Sigma$ is continuous it
follows that $\varphi(x_{i})\rightarrow0=\varphi(x)$, when $x_{i}\rightarrow x\in\partial M$. 

Next, fix a compact subset $K$ in the interior of $M$ and $\varepsilon>0.$
Let \hbox{$M_{\delta}:=\{\eta_\Sigma<-\delta\}$} where $\delta$ is sufficiently
small to make sure that $K$ is contained in $M_{4\delta}$. By the
regularization results in [Dem92] or [Dem94], there is a sequence $\varphi_j$
in $\PSH(M,\alpha-2^{-j}\alpha)\cap C^0(M_{\delta/2})$ decreasing
to the upper semi-continuous regularization~$\varphi^*$. By replacing
$\varphi_j$ with $(1-2^{-j})^{-1}\varphi_j$, we can even assume
$\varphi_j\in \PSH(M,\alpha)\cap C^0(M_{\delta/2})$. Put 
$$
\varphi_j':=\max\{\varphi_j-\varepsilon,C\eta_\Sigma\}~~
\hbox{on $M_{\delta}$},\quad\hbox{and}\quad 
\varphi_j':=C\eta_\Sigma~~\hbox{on $M\smallsetminus M_{\delta}$.}
$$
On $\partial M_\delta$ we have $C\eta_\Sigma=-C\delta$ and we can take $j$ 
so large that 
$$
\varphi_j<-C\eta_\Sigma+\varepsilon/2=C\delta+\varepsilon/2,
$$
so we will have $\varphi_j-\varepsilon<C\eta_\Sigma$ as soon as 
$2C\delta\leq\varepsilon/2$. We simply take $\varepsilon=4C\delta$.
Then $\varphi_j'$ is a well defined continuous $\alpha$-psh function on $M$, and
$\varphi'_j$ is equal to $\varphi_j-\varepsilon$ on $K\subset M_{4\delta}$
as $C\eta_\Sigma\le -4C\delta\le-\varepsilon\leq\varphi_j-\varepsilon$ there.
In particular, $\varphi_j'$ is a candidate
for the sup defining $\varphi$, hence
$\varphi'_j\leq\varphi\leq\varphi^*$ and so
$$
\varphi^*\leq\varphi_j\leq\varphi'_j+\varepsilon\leq\varphi^*+\varepsilon
$$
on $K$. This means that $\varphi_j$ converges to $\varphi$ uniformly
on $K$ and therefore $\varphi$ is continuous on $K$. All in all this shows
that $\varphi\in C^0(M)$. The last statement of the proposition
follows from standard local considerations for envelopes due to 
Bedford-Taylor [BT76] (see also the exposition made in [Dem89]).\qed

\claim (4.4) Theorem|Let $\alpha$ be a closed real $(1,1)$-form on $M$ 
with bounded  coefficients, such that $\alpha_{|\{s\}\times X}
\ge\varepsilon_0\omega$ is positive definite for all $s\in\Sigma$.
Consider a continuous function $f$ on $\partial M$
such that $f_s\in\PSH(X,\alpha_s)$
for all $s\in\partial\Sigma$. Then the upper envelope 
$\varphi=\varphi_{\alpha,f}$ 
is the unique $\alpha$-psh continuous solution of the Dirichlet problem
$$
\varphi=f~~\hbox{on}~\partial M,\qquad
(dd^cu+\alpha)^{\dim M}=0\,\,\,\hbox{on the interior}~M^\circ.\leqno(4.5)
$$
 Moreover, if $f$ is $C^{1,1}$ on $\partial M$ then,
for any $s$ in $\Sigma$, the restriction $\varphi_s$ of $\varphi$ 
on $\{s\}\times X$
has a $dd^c$ in $L_{\loc}^{\infty}$. More precisely, we have
a uniform bound $|dd^c\varphi_s|_\omega\leq C$
a.e.\ on $X$, where $C$ is a constant independent of $s$.
\endclaim

\proof. Without loss of generality, we may assume as in Lemma (4.4)
that $\alpha$ is positive definite on~$M$. Also, after adding a positive 
constant to~$f$, which only has the effect of adding the same constant
to $\varphi=\varphi_{\alpha,f}$, we may suppose that 
$\sup_{\partial M} f>0$ (this will simplify a little bit the 
arguments below).

\medskip\noindent
\emph{Continuity.} Let us first prove the continuity statement in
the theorem. In the case when $f$ extends to a smooth function $F$ in
$\PSH(M,(1-\varepsilon)\alpha)$ the statement follows immediately from 
(4.2) and Lemma~(4.3) since
$$
f-F=0~~\hbox{on $\partial M$ and}~~\beta=\alpha+dd^c F\ge
\varepsilon\alpha\ge\varepsilon\varepsilon_0\omega.
$$
Next, assume that $f$ is \emph{smooth} on $\partial M$
and that $f_s\in\PSH(X,(1-\varepsilon)\alpha_s)$ for all 
$s\in\partial\Sigma$. If we take a smooth extension 
$\widetilde f$ of $f$ to $M$ and $C\gg 1$, we will get
$$
\alpha+dd^c(\widetilde f(x,s)+C\eta_\Sigma(s))\ge (\varepsilon/2)\alpha
$$
on a sufficiently small neighborhood $V$ of $\partial M$ (again by
using Cauchy-Schwarz). Therefore, after enlarging $C$ if necessary, 
we can define
$$
F(x,s)=\max\nolimits_\varepsilon(\widetilde f(x,s)+C\eta_\Sigma(s),0)
$$
with a regularized max function $\max\nolimits_\varepsilon$, in
such a way that the maximum is equal to $0$ on a neighborhood of $M\ssm V$
($C\gg 1$ being used to ensure that $\widetilde f+C\eta_\Sigma<0$ on $M\ssm V$).
Then $F$ equals $f$ on $\partial M$ and satisfies 
$$
\alpha+dd^cF\ge (\varepsilon/2)\alpha\ge(\varepsilon\varepsilon_0/2)\omega
$$
on $M$, and we can argue as previously. Finally, to handle the general 
case where $f$ is continuous with $f_s\in\PSH(X,\alpha_s)$ for every
$s\in\Sigma$, we may, by a parametrized version of
Richberg's regularization theorem applied to 
$(1-2^{-\nu})f+C\,2^{-\nu}$ (see e.g.\ [Dem91]), write $f$ as 
a decreasing uniform limit of smooth
functions $f_{\nu}$ on $\partial M$ satisfying 
$f_{\nu,s}\in\PSH(X,(1-2^{-\nu-1})\alpha_s)$ for every $s\in
\partial\Sigma$. Then
$\varphi_{\omega,f}$ is a decreasing uniform limit on $M$ of the 
continuous functions $\varphi_{\omega,f_{\nu}}$
(as follows easily from the definition of $\varphi_{\omega,f}$ as 
an upper envelope). Observe also that the uniqueness of
a continuous solution of the Dirichlet problem (4.5)
results from a standard application of the maximum principle for the
Monge-Amp\`ere operator. This proves the general case of
the continuity statement.

\medskip\noindent
\emph{Smoothness.} Next, we turn to the proof of the smoothness statement.
Since the proof is a straightforward adaptation of the proof of the
main regularity result above we will just briefly indicate the relevant
modification. Quite similarly to what we did in section~1, we consider 
an $\alpha$-psh function
$\psi$ with $\psi\le f$ on $\partial M$, and introduce  the fiberwise
transform $\Psi_s$ of $\psi_s$ on each $\{s\}\times X$ which is
defined in terms of the exponential map~$\exph:T_X\rightarrow X$, and
we put
$$
\Psi(z,s,t)=\Psi_s(z,t).
$$
 Then essentially the same calculations as in the previous case show
that all properties of $\Psi$ are still valid with the constant $K$
depending on the $C^{1,1}$-norm of the local potentials $u(z,s)$
of $\alpha$, the constant $A$ depending only on $\omega$ and
with
$$
\partial\Psi(z,s,t)/\partial(\log t):=\lambda(z,s,t)\rightarrow\nu(\psi_s),
$$
 as $t\rightarrow0^+$, where $\nu(\psi_s)$ is the Lelong number
of the function $\psi_s$ on $X$ at~$z$. Moreover, the local vector
valued differential $dz$ should be replaced by the differential $d(z,s)=dz+ds$
in the previous formulas. Next, performing a Kiselman-Legendre transform
fiberwise we let
$$
\psi_{c,\delta}(z,s):=(\psi_s)_{c,\delta}(z)
$$
Then, using a parametrized version of the estimates of [Dem94] and 
the properties of $\Psi(z,s,w)$ as in section~1, arguments
derived from Kiselman's infimum principle show that
$$
\alpha+dd^c\psi_{c,\delta}\geq(-A\min(c,\lambda(z,s,\delta))-K\delta^2)\omega_M
\geq-(Ac+K\delta^2)\omega_M,
\leqno(4.6)
$$
where $\omega_M$ is the K\"ahler form on $M$. In addition to this, we
have $|\psi_{c,\delta}-f|\le K'\delta^2$ on $\partial M$ by the hypothesis that
$f$ is $C^{1,1}$. For a sufficiently large constant $C_1$, we infer from this
that $\theta=(1-C_1(Ac+K\delta^2))\psi_{c,\delta}$ satisfies
$\theta\le f$ on $\partial M$ (here we use the fact that $f>0$ and hence
that $\psi_0\equiv 0$ is a candidate for the upper envelope). Moreover
$\alpha+dd^c\theta\ge 0$ on $M$ thanks to (4.6) and the positivity 
of~$\alpha$. Therefore $\theta$ is a candidate for the upper envelope and so
$\theta\le \varphi=\varphi_{f,\alpha}$. Repeating the arguments of section~1
almost word by word, we obtain for $(\rho_t\varphi)(z,s):=\Phi(z,s,t)$
the analogue of estimate (1.15) which reduces simply to
$$
\liminf_{t\to 0_+}{\rho_t\varphi(z,s)-\varphi(z,s)\over t^2}\leq C_2,
$$
as $\psi_0\equiv 0$ in the present situation. The final conclusion
follows from (1.16) and the related arguments already explained.\qed
\medskip

In connection to the study of Wess-Zumino-Witten type equations 
[Don99], [Don02]
and geodesics in the space of K\"ahler metrics [Don99], [Don02], [Che00] it
is useful to formulate the result of the previous theorem as an extension
problem from $\partial\Sigma$, in the case when $\alpha(z,s)=\omega(z)$ does
not depend on~$s$. To this end, 
let $F:\,\partial\Sigma\to\PSH(X,\omega)$
be the map defined by $F(s)=f_s$. Then the previous theorem gives
a continuous ``maximal plurisubharmonic'' extension $U$ of $F$
to $\Sigma$, where $U(s):=u_s$ so that
$U:\,\partial\Sigma\to\PSH(X,\omega)$. 

Let us next specialize to the case when $\Sigma:=A$ is an annulus
$R_1<|s|<R_2$ in~$\bC$ and the boundary data $f(x,s)$ is invariant 
under rotations $s\mapsto s\,e^{i\theta}$.
Denote by $f^0$ and $f^1$ the elements in $\PSH(X,\omega)$
corresponding to the two boundary circles of~$A$. Then the previous
theorem furnishes a continuous path $f^t$ in $\PSH(X,\omega)$,
if we put $t=\log|s|$, or rather $t=\log(|s|/R_1)/\log(R_2/R_1)$
to be precise. Following [PS08] the corresponding path of semi-positive 
forms $\omega^t:=\omega+dd^cf^t$
will be called a \emph{$($generalized$)$ geodesic} in $\PSH(X,\omega)$
(compare also with Remark 4.8).

\claim (4.7) Corollary|
Assume that the semi-positive closed $(1,1)$ forms $\omega^0$
and $\omega^1$ belong to the same K\"ahler class $\{\omega\}$ and
have bounded coefficients. Then the geodesic $\omega^t$
connecting $\omega^0$ and $\omega^1$ is continuous on $[0,1]\times X$,
and there is a constant $C$ such that $\omega^t\leq C\omega$ on
$X$, i.e.\ $\omega^t$ has uniformly bounded coefficients.
\endclaim

In particular, the previous corollary shows that the space of all
semi-positive forms with bounded coefficients, in a given K\"ahler class,
is ``geodesically convex''.

\claim (4.8) Remark|{\rm 
As shown in the work of Semmes, Mabuchi and Donaldson, the space of
K\"ahler metrics $\cH_\omega$ in a given K\"ahler class $\{\omega\}$
admits a natural Riemannian structure defined in the following way
(see [Che00] and references therein). First note that
the map $u\mapsto\omega+dd^cu$ identifies $\cH_\omega$
with the space of all smooth and strictly $\omega$-psh functions,
modulo constants. Now identifying the tangent space of $\cH_\omega$
at the point $\omega+dd^cu\in\cH_\omega$ with $C^\infty(X)/\bR$,
the squared norm of a tangent vector $v$ at the point $u$ is defined
as 
$$
\int_{X}v^{2}(\omega+dd^cu)^n/n!.
$$
Then the potentials $f^t$ of any given geodesic $\omega^t$ in
$\cH_\omega$ are in fact solutions of the Dirichlet problem
(4.5) above, with $\Sigma$ an annulus and $t:=\log|s|$, see [Che00]. 
However, the \emph{existence} of a geodesic $u_t$ in
$\cH_\omega$ connecting any given points $u_0$ and $u_1$
is an open and even dubious problem. In the case when $\Sigma$ is
a Riemann surface, the boundary data $f$ is smooth with 
$\alpha_s+dd^cf_s>0$
on $X$ for $s\in\partial\Sigma$ it was shown in [Che00]
that the solution $\varphi$ of the Dirichlet problem (4.5)
has a total Laplacian which is bounded on~$M$.
See also [Blo08] for a detailed analysis of the proof in [Che00]
and some refinements. On the other hand it is not known
whether $\alpha_s+dd^c\varphi_s>0$ for all $s\in\Sigma$, even
under the assumption of rotational invariance which appears in the
case of geodesics as above. However, see [CT08] for results in
this direction. A~case similar to the degenerate setting in the previous
corollary was also considered very recently in [PS08], building
on [Blo08].}
\endclaim

\claim (4.9) Remark|{\rm 
Note that the assumption $f\in C^2(\partial M)$ is not sufficient
to obtain uniform estimates on the total Laplacian on $M$ with respect 
to $\omega_M$
of the envelope $u$ up to the boundary. To see this let $\Sigma$
be the unit-ball in $\bC^2$ and write $s=(s_1,s_2)\in\bC^2$.
Then $f(s):=(1+\Re s_1)^{2-\varepsilon}$ is in $C^{4-2\varepsilon}(\partial M)$
and $u(x,s):=f(s)$ is the continuous solution of the Dirichlet problem
(4.5). However, $u$ is not in $C^{1,1}(M)$
at $(x\,;\,-1,0)\in\partial M$ for any $x\in X$. Note that this exemple
is the trivial extension of the exemple in [CNS86] for the real
Monge-Amp\`ere equation in the disc.}
\endclaim

\section{5. Regularity of ``supercanonical'' metrics}

Let $X$ be a compact complex manifold and $(L,h_{L,\gamma})$ a holomorphic
line bundle over $X$ equipped with a singular hermitian metric
$h_{L,\gamma}=e^{-\gamma}h_L$ with satisfies $\int e^{-\gamma}<+\infty$
locally on~$X$, where $h_L$ is a smooth metric on~$L$. In fact, we can 
more generally consider the case where
$(L,h_{L,\gamma})$ is a ``hermitian $\bR$-line bundle''; by this we mean that
we have chosen a smooth real $d$-closed $(1,1)$ form $\alpha_L$ on~$X$
(whose $dd^c$ cohomology class is equal to $c_1(L))$, and a specific 
current $T_{L,\gamma}$ representing it, namely
$T_{L,\gamma}=\alpha_L+dd^c\gamma$, such that $\gamma$ is a locally
integrable function satisfying~$\int e^{-\gamma}<+\infty$. An important
special case is obtained by considering a klt (Kawamata log terminal) 
effective divisor~$\Delta$. In this situation
$\Delta=\sum c_j\Delta_j$ with $c_j\in\bR$, and if $g_j$ is a local
generator of the ideal sheaf $\cO(-\Delta_j)$ identifying it to the trivial
invertible sheaf $g_j\cO$, we take
$\gamma=\sum c_j\log|g_j|^2$, $T_{L,\gamma}=\sum c_j[\Delta_j]$
(current of integration on $\Delta$) and $\alpha_L$ given by any smooth
representative of the same $dd^c$-cohomology class; the klt condition
precisely means that
$$
\int_V e^{-\gamma}=\int_V \prod|g_j|^{-2c_j}<+\infty
\leqno(5.1)
$$ 
on a small neighborhood~$V$ of any point in the support 
$|\Delta|=\bigcup\Delta_j$ (condition (5.1) implies $c_j<1$
for every $j$, and this in turn is sufficient to imply $\Delta$ klt
if $\Delta$ is a normal crossing divisor; the line bundle $L$ is 
then the real line bundle $\cO(\Delta)$, which makes sens as a
genuine line bundle only if $c_j\in\bZ$).
For~each klt pair $(X,\Delta)$ such that $K_X+\Delta$
is pseudo-effective, H.~Tsuji [Ts07a, Ts07b] has introduced a
``supercanonical metric'' which generalizes the metric introduced
by Narasimhan and Simha [NS68] for projective algebraic varieties
with ample canonical divisor. We take the opportunity to present
here a simpler, more direct and more general approach. 

We assume from now on that $K_X+L$ is \emph{pseudo-effective}, i.e.\ 
that the class
$c_1(K_X)+\{\alpha_L\}$ is pseudo-effective, and under this condition, we
are going to define a ``supercanonical metric'' on $K_X+L$. Select
an arbitrary smooth hermitian metric $\omega$ on $X$. We then find 
induced hermitian
metrics $h_{K_X}$ on $K_X$ and $h_{K_X+L}=h_{K_X}h_L$ on $K_X+L$, whose 
curvature is the smooth real $(1,1)$-form
$$
\alpha=\Theta_{K_X+L,h_{K_X+L}}=\Theta_{K_X,\omega}+\alpha_L.
$$
A singular hermitian metric
on $K_X+L$ is a metric of the form $h_{K_X+L,\varphi}=e^{-\varphi}h_{K_X+L}$
where $\varphi$ is locally
integrable, and by the pseudo-effectivity assumption, we can find 
quasi-psh functions $\varphi$ such that $\alpha+dd^c\varphi\ge 0$.
The metrics on $L$ and $K_X+L$ can now be ``subtracted'' to give rise
to a metric
$$
h_{L,\gamma}h_{K_X+L,\varphi}^{-1}=e^{\varphi-\gamma}h_Lh_{K_X+L}^{-1}
=e^{\varphi-\gamma}h_{K_X}^{-1}=e^{\varphi-\gamma}dV_\omega
$$
on $K_X^{-1}=\Lambda^nT_X$, since $\smash{h_{K_X}^{-1}}=dV_\omega$ is just
the hermitian $(n,n)$ volume form on~$X$. Therefore the integral
$\int_X h_{L,\gamma}h_{K_X+L,\varphi}^{-1}$ has an intrinsic meaning, and 
it makes sense to require that
$$
\int_X h_{L,\gamma}h_{K_X+L,\varphi}^{-1}
=\int_X e^{\varphi-\gamma}dV_\omega\le 1\leqno(5.2)
$$
in view of the fact that $\varphi$ is locally bounded from above and of
the assumption $\int e^{-\gamma}<+\infty$. Observe that condition (5.2) can
always be achieved by subtracting a constant to $\varphi$. Now, 
we can generalize Tsuji's supercanonical metrics on klt pairs
(cf.\ [Ts07b]) as follows.

\claim (5.3) Definition|Let $X$ be a compact 
complex manifold and let $(L,h_L)$ be a hermitian $\bR$-line bundle on $X$
associated with a smooth real closed $(1,1)$ form~$\alpha_L$.
Assume that $K_X+L$ is pseudo-effective and that $L$ is
equipped with a singular hermitian metric $h_{L,\gamma}=e^{-\gamma}h_L$ 
such that \hbox{$\int e^{-\gamma}<+\infty$} locally on~$X$. Take a hermitian
metric $\omega$ on $X$ and define 
$\alpha=\Theta_{K_X+L,h_{K_X+L}}=\Theta_{K_X,\omega}+\alpha_L$.
Then we define the supercanonical metric $h_\can$ of $K_X+L$ to be 
$$
\eqalign{
&h_{K_X+L,\can}=\inf_\varphi h_{K_X+L,\varphi}\quad
\hbox{i.e.}\quad h_{K_X+L,\can}=e^{-\varphi_{\can}}h_{K_X+L},~~
\hbox{where}\cr
&\varphi_{\can}(x)=\sup_\varphi\varphi(x)~~
\hbox{for all $\varphi$ with}~~\alpha+dd^c\varphi\ge 0,
~~\int_X e^{\varphi-\gamma}dV_\omega\le 1.\cr}
$$
\endclaim

In particular, this gives a definition of the supercanonical metric
on $K_X+\Delta$ for every klt pair $(X,\Delta)$ such that
$K_X+\Delta$ is pseudo-effective, and as an even more special case,
a supercanonical metric on $K_X$ when $K_X$ is pseudo-effective.

In the sequel, we assume that $\gamma$ has analytic singularities, otherwise
not much can be said. The mean value inequality then immediately shows that the 
quasi-psh functions $\varphi$ involved in definition (5.3) 
are globally uniformly bounded outside of the poles of $\gamma$, and
therefore everywhere on $X$, hence the envelopes
$\varphi_\can=\sup_\varphi\varphi$ are indeed well defined and bounded
above. As a consequence, we get a ``supercanonical'' current
$T_\can=\alpha+dd^c\varphi_\can\ge 0$ and $h_{K_X+L,\can}$ satisfies
$$
\int_Xh_{L,\gamma}h_{K_X+L,\can}^{-1}=
\int_X e^{\varphi_\can-\gamma}dV_\omega<+\infty.\leqno(5.4)
$$
It is easy to see that in Definition (5.3) the supremum is a maximum and that 
$\varphi_\can=(\varphi_\can)^*$ every\-where, so that
taking the upper semicontinuous regularization is not needed. In fact
if $x_0\in X$ is given and we write
$$
(\varphi_\can)^*(x_0)=\limsup_{x\to x_0}\varphi_\can(x)
=\lim_{\nu\to+\infty}\varphi_\can(x_\nu)
=\lim_{\nu\to+\infty}\varphi_\nu(x_\nu)
$$
with suitable sequences $x_\nu\to x_0$ and $(\varphi_\nu)$ such that
$\int_Xe^{\varphi_\nu-\gamma}dV_\omega\le 1$, the well-known 
weak compactness properties of quasi-psh functions in  $L^1$ topology 
imply the existence of
a subsequence of $(\varphi_\nu)$ converging in $L^1$ and almost everywhere to 
a quasi-psh limit $\varphi$. Since $\int_Xe^{\varphi_\nu-\gamma}dV_\omega
\le 1$ holds true for every~$\nu$, Fatou's lemma 
implies that we have $\int_Xe^{\varphi-\gamma}dV_\omega\le 1$ in
the limit. By taking a subsequence, we can 
assume that $\varphi_\nu\to\varphi$ in $L^1(X)$. Then for every $\varepsilon>0$
the mean value $-\kern-8.7pt\int_{B(x_\nu,\varepsilon)}\varphi_\nu$ 
satisfies
$$
-\kern-10.7pt\int_{B(x_0,\varepsilon)}\varphi=\lim_{\nu\to +\infty}
-\kern-10.7pt\int_{B(x_\nu,\varepsilon)}\varphi_\nu\ge \lim_{\nu\to+\infty}
\varphi_\nu(x_\nu)=(\varphi_\can)^*(x_0),
$$
hence we get $\varphi(x_0)=\lim_{\varepsilon\to 0}
-\kern-8.7pt\int_{B(x_0,\varepsilon)}\varphi\ge 
(\varphi_\can)^*(x_0)\ge \varphi_\can(x_0)$, and therefore the 
sup is a maximum and $\varphi_\can=\varphi_\can^*$.

By elaborating on this argument, we can infer certain regularity
properties of the envelope. However, there is no reason why the
integral occurring in (5.4) should be equal to~$1$
when we take the upper envelope. As a consequence, neither the 
upper envelope nor its regularizations participate to the family
of admissible metrics. This is the reason why the estimates that we
will be able to obtain are much weaker than in the case of envelopes
normalized by a condition~$\varphi\le 0$.

\claim (5.5) Theorem|Let $X$ be a compact complex manifold and $(L,h_L)$ a 
holomorphic $\bR$-line bundle such that $K_X+L$ is big. Assume that
$L$ is equipped with a singular hermitian metric $h_{L,\gamma}=e^{-\gamma}h_L$
with analytic singularities such that $\int e^{-\gamma}<+\infty$
$($klt condition$)$. Denote by $Z_0$ the set of poles of a singular 
metric $h_0=e^ {-\psi_0}h_{K_X+L}$ with analytic singularities 
on~$K_X+L$ and by $Z_\gamma$ the poles of~$\gamma$ $($assumed analytic$)$.
Then the associated supercanonical metric $h_\can$ is continuous on 
$X\ssm (Z_0\cup Z_\gamma)$ and possesses some computable logarithmic 
modulus of continuity.
\endclaim

\proof. With the notation already introduced, let 
$h_{K_X+L,\varphi}=e^{-\varphi}h_{K_X+L}$ be a singular hermitian metric
such that its curvature satisfies $\alpha+dd^c\varphi\ge 0$ and 
$\int_Xe^{\varphi-\gamma}dV_\omega\le 1$. We apply to~$\varphi$ the 
regularization procedure defined in (1.6). Jensen's inequality implies
$$
e^{\Phi(z,w)}\le 
\int_{\zeta\in T_{X,z}}e^{\varphi(\exph_z(w\zeta))}\,\chi(|\zeta|^2)\,
dV_\omega(\zeta).
$$
If we change variables by putting $u=\exph_z(w\zeta)$,
then in a neighborhood of the diagonal of $X\times X$ we have an inverse map
$\logh:X\times X\to T_X$ such that $\exph_z(\logh(z,u))=u$ and we find
for $w$ small enough
$$
\eqalign{
\int_X &e^{\Phi(z,w)-\gamma(z)}dV_\omega(z)\cr
&\le\int_{z\in X}\bigg(\int_{u\in X}
e^{\varphi(u)-\gamma(z)}\chi\bigg({|\logh(z,u)|^2\over |w|^2}\bigg)
{1\over |w|^{2n}}dV_\omega(\logh(z,u))\bigg)dV_\omega(z)\cr
&=\int_{u\in X}P(u,w)\,e^{\varphi(u)-\gamma(u)}dV_\omega(u)\cr}
$$
where  $P$ is a kernel on $X\times D(0,\delta_0)$ such that
$$
P(u,w)=
\int_{z\in X} {1\over |w|^{2n}}
\chi\bigg({|\logh(z,u)|^2\over |w|^2}\bigg)
{e^{\gamma(u)-\gamma(z)}dV_\omega(\logh(z,u))\over dV_\omega(u)}\,dV_\omega(z).
$$
Let us first assume that $\gamma$ is smooth (the case where $\gamma$ has
logarithmic poles will be considered later). Then a change of variable
$\zeta={1\over w}\logh(z,u)$ shows that $P$ is smooth
and we have $P(u,0)=1$. Since $P(u,w)$ depends only on $|w|$ we infer
$$
P(u,w)\le 1+C_0|w|^2
$$ 
for $w$ small. This shows that the integral
of $z\mapsto e^{\Phi(z,w)-C_0|w|^2}$ will be at most equal to $1$,
and therefore if we define 
$$
\varphi_{c,\delta}(z)=\inf_{t\in{}]0,\delta]}
\Phi(z,t)+Kt^2-K\delta^2-c\log{t\over\delta}\leqno(5.6)
$$
as in (1.10), the function
$\varphi_{c,\delta}(z)\le \Phi(z,\delta) $ will also satisfy
$$
\int_Xe^{\varphi_{c,\delta}(z)-C_0\delta^2-\gamma(z)}dV_\omega\le 1.\leqno(5.7)
$$
Now, thanks to the assumption that $K_X+L$ is big, there exists a quasi-psh
function $\psi_0$ with analytic singularities such that
$\alpha+dd^c\psi_0\ge \varepsilon_0\omega$. We can assume
$\int_Xe^{\psi_0-\gamma}dV_\omega=1$ after adjusting $\psi_0$ with
a suitable constant. Consider a pair of points $x,y\in X$. We 
take $\varphi$ so that $\varphi(x)=\varphi_\can(x)$ (this is possible by the 
above discussion). We define
$$
\varphi_\lambda=\log\big(\lambda e^{\psi_0}+(1-\lambda)e^\varphi\big)\leqno(5.8)
$$
with a suitable constant $\lambda\in[0,1/2]$ which will be fixed later,
and obtain in this way regularized functions $\Phi_\lambda(z,w)$
and $\varphi_{\lambda,c,\delta}(z)$. This is obviously a compact family and
therefore the associated constants $K$ needed in (5.6) are uniform 
in~$\lambda$. Also, as in section~1, we have
$$
\alpha+dd^c\varphi_{\lambda,c,\delta}\ge -(Ac+K\delta^2)\,\omega\qquad
\hbox{for all $\delta\in{}]0,\delta_0]$}.\leqno(5.9)
$$
Finally, we consider the linear combination
$$
\theta={Ac+K\delta^2\over\varepsilon_0}\psi_0+
\bigg(1-{Ac+K\delta^2\over\varepsilon_0}\bigg)
(\varphi_{\lambda,c,\delta}-C_0\delta^2).\leqno(5.10)
$$
Clearly, $\int_X e^{\varphi_\lambda-\gamma}dV_\omega\le 1$, and therefore $\theta$
also satisfies $\int_Xe^{\theta-\gamma}dV_\omega\le 1$ by H\"older's inequality.
Our linear combination is precisely taken so that $\alpha+dd^c\theta\ge 0$.
Therefore, by definition of $\varphi_\can$, we find that
$$
\varphi_\can\ge\theta={Ac+K\delta^2\over\varepsilon_0}\psi_0+
\bigg(1-{Ac+K\delta^2\over\varepsilon_0}\bigg)
(\varphi_{\lambda,c,\delta}-C_0\delta^2).
\leqno(5.11)
$$
Assume $x\in X\ssm Z_0$, so that $\varphi_\lambda(x)>-\infty$ and
$\nu(\varphi_\lambda,x)=0$. In (5.6), the infimum is reached either 
for $t=\delta$ or for $t$ such that $c=t{d\over dt}(\Phi_\lambda(z,t)+Kt^2)$.
The function  $t\mapsto \Phi_\lambda(z,t)+Kt^2$ is convex increasing
in $\log t$ and tends to $\varphi_\lambda(z)$ as $t\to 0$. By~convexity,
this implies
$$
\eqalign{
c=t{d\over dt}(\Phi_\lambda(z,t)+Kt^2)
&\le {(\Phi_\lambda(x,\delta_0)+K\delta_0^2)
-(\Phi_\lambda(z,t)+Kt^2)\over \log(\delta_0/t)}\cr
&\le{C_1-\varphi_\lambda(x)\over \log(\delta_0/t)}
\le{C_1+|\psi_0(z)|+\log(1/\lambda)\over \log(\delta_0/t)},\cr}
$$
hence
$$
{1\over t}\le \max\bigg({1\over\delta}\,,\;{1\over\delta_0}
\exp\Big({C_1+|\psi_0(z)|+\log(1/\lambda)\over c}\Big)\bigg).
\leqno(5.12)
$$
This shows that $t$ cannot be too small when the infimum is reached. When
$t$ is taken equal to the value which achieves the infimum for 
$z=y$, we find
$$
\varphi_{\lambda,c,\delta}(y)=\Phi_\lambda(y,t)+Kt^2-K\delta^2-c\log{t\over\delta}
\ge\Phi_\lambda(y,t)+K t^2-K\delta^2.
\leqno(5.13)
$$
Since $z\mapsto\Phi_\lambda(z,t)$ is a convolution of $\varphi_\lambda$,
we get a bound of the first order derivative
$$
|D_z\Phi_\lambda(z,t)|
\le \Vert\varphi_\lambda\Vert_{L^1(X)} {C_2\over t} \le {C_3\over t},
$$
and with respect to the geodesic distance $d(x,y)$ we infer from this
$$
\Phi_\lambda(y,t)\ge \Phi_\lambda(x,t)-{C_3\over t}d(x,y).
\leqno(5.14)
$$
A combination of (5.11), (5.13) and (5.14) yields
$$
\eqalign{
\varphi_\can(y)
&\ge {Ac+K\delta^2\over\varepsilon_0}\psi_0(y)+
\bigg(1-{Ac+K\delta^2\over\varepsilon_0}\bigg)
\Big(\Phi_\lambda(x,t){+}Kt^2{-}K\delta^2{-}{C_3\over t}d(x,y)\Big)\cr
&\ge {Ac+K\delta^2\over\varepsilon_0}\psi_0(y)+\bigg(1-{Ac+K\delta^2\over
\varepsilon_0}\bigg)\Big(\varphi_\lambda(x)-K\delta^2-{C_3\over t}d(x,y)\Big)\cr
&\ge \log\big(\lambda e^{\psi_0(x)}+(1-\lambda)e^{\varphi(x)}\big)-
C_4\Big((c+\delta^2)(|\psi_0(y)|+1)+{1\over t}d(x,y)\Big),\cr
&\ge \varphi_\can(x)-
C_5\Big(\lambda+(c+\delta^2)(|\psi_0(y)|+1)+{1\over t}d(x,y)\Big),\cr
}
$$
if we use the fact that $\varphi_\lambda(x)\le C_6$, 
$\varphi(x)=\varphi_\can(x)$ and
$\log(1-\lambda)\ge-(2\log 2)\lambda$ for all $\lambda\in[0,1/2]$. 
By exchanging the roles of $x,y$ and using (5.12), we see that for all $c>0$,
$\delta\in{}]0,\delta_0]$ and $\lambda\in{}]0,1/2]$, there is an inequality
$$
\big|\varphi_\can(y)-\varphi_\can(x)\big|\le C_5
\Big(\lambda+(c+\delta^2)\big(\max(|\psi_0(x)|,|\psi_0(y)|)+1\big)
+{1\over t}d(x,y)\Big)\leqno(5.15)
$$
where
$$
{1\over t}\le \max\bigg({1\over\delta}\,,\;{1\over\delta_0}
\exp\Big({C_1+\max(|\psi_0(x)|,|\psi_0(y)|)+\log(1/\lambda)\over c}\Big)\bigg).
\leqno(5.16)
$$
By taking $c$, $\delta$ and $\lambda$ small, one easily sees that this
implies the continuity of $\varphi_\can$ on $X\ssm Z_0$. More precisely,
if we choose
$$
\delta=d(x,y)^{1/2},~\lambda={1\over |\log d(x,y)|},~
c={C_1{+}\max(|\psi_0(x)|,|\psi_0(y)|){+}
\big|\log|\log d(x,y)|\big|\over \log \delta_0/d(x,y)^{1/2}}\!
$$
with $d(x,y)<\delta_0^2<1$, we get ${1\over t}\le d(x,y)^{-1/2}$, whence
an explicit, but certainly non optimal, modulus of continuity of the form
$$
\big|\varphi_\can(y)-\varphi_\can(x)\big|\le 
C_7\big(\max(|\psi_0(x)|,|\psi_0(y)|)+1\big)^2
{\big|\log|\log d(x,y)|\big|+1\over |\log d(x,y)|+1}\;.
$$
When the weight $\gamma$ has analytic singularities, the kernel $P(u,w)$ is no
longer smooth and the volume estimate (5.7). In this case, we
use a modification $\mu:\smash{\wh X}\to X$ in such a way that the singularities
of $\gamma\circ\mu$ are divisorial, given by a divisor with normal
crossings. If we put 
$$
\wh L=\mu^*L-K_{\wh X/X}=\mu^* L-E
$$ 
($E$ the exceptional divisor), then we get an induced singular metric 
on $\wh L$ which still satisfies the klt condition, and the corresponding
supercanonical metric on $K_{\wh X}+\wh L$ is just the pull-back by
$\mu$ of the supercanonical metric on $K_X+L$. This shows that we may assume
from the start that the singularities of $\gamma$ are divisorial and
given by a klt divisor~$\Delta$. In this
case, a solution to the problem is to introduce
a complete hermitian metric $\hat\omega$ of uniformly bounded curvature on
$X\ssm|\Delta|$ by using the Poincar\'e metric on the punctured disc as a 
local model transversally to the components of~$\Delta$. The Poincar\'e
metric on the punctured unit disc is given by
$$
{|dz|^2\over |z|^2(\log|z|)^2}
$$
and the singularity of $\hat\omega$ along the component 
$\Delta_j=\{g_j(z)=0\}$ of $\Delta$ is given by 
$$
\hat\omega=\sum -dd^c\log|\log|g_j||\quad\mod C^\infty. 
$$
Since such
a metric has bounded geometry and this is all that we need for the
calculations of [Dem94] to work, the estimates that we have made here
are still valid,
especially the crucial lower bound $\alpha+dd^c\varphi_{\lambda,c,\delta}\ge
-(Ac+K\delta^2)\,\hat\omega$. In order to compensate this loss of positivity, we need
a quasi-psh function $\smash{\hat\psi_0}$ such that
$\alpha+dd^c\hat\psi_0\ge \varepsilon_0\hat\omega$, but such a lower
bound is possible by adding terms of the form $-\varepsilon_1\log|\log|g_j||$
to our previous quasi-psh function $\psi_0$. Now, with respect to the
Poincar\'e metric, a $\delta$-ball of center $z_0$ in the punctured disc
is contained in the corona 
$$
|z_0|^{e^{-\delta}}<|z|<|z_0|^{e^\delta},
$$
and it is easy to see from there that the mean value of $|z|^{-2a}$ on
a $\delta$-ball of center $z_0$ is multiplied by at most 
$|z_0|^{-2a\delta}$. This implies that a function of
the form 
$\hat\varphi_{c,\delta}=\varphi_{c,\delta}+C_9\delta\sum\log|g_j|$ will 
actually give 
rise to an integral $\int_Xe^{\hat\varphi_{c,\delta}-\gamma}dV_\omega\le 1$.
We see that the term $\delta^2$ in (5.15) has to
be replaced by a term of the form 
$$
\delta\sum\max\big(|\log|g_j(x)||,|\log|g_j(x)||\big).
$$
This is enough to obtain the continuity of $\varphi_\can$ on
$X\ssm(Z_0\cup|\Delta|)$, as well as an explicit logarithmic modulus
of continuity.\qed

\claim (5.17) Algebraic version|{\rm Since the klt condition is open 
and $K_X+L$ is assumed to be~big, we can always perturb $L$ a 
little bit, and after
blowing-up $X$, assume that $X$ is projective and that
$(L,h_{L,\gamma})$ is obtained as a sum of $\bQ$-divisors
$$
L=G+\Delta
$$
where $\Delta$ is klt and $G$ is equipped with a smooth
metric $h_G$ (from which $h_{L,\gamma}$ is inferred, with $\Delta$
as its poles, so that $\Theta_{L,h_{L,\gamma}}=\Theta_{G,L_G}+[\Delta]$).
Clearly this situation is ``dense'' in what we
have been considering before, just as $\bQ$ is dense in $\bR$.
In this case, it is possible to 
give a more algebraic definition of the supercanonical
metric $\varphi_\can$, following the original idea of Narasimhan-Simha
[NS68] (see also H.~Tsuji [Ts07a])~--~the case considered by these
authors is the special situation where $G=0$, $h_G=1$ (and moreover
$\Delta=0$ and $K_X$ ample, for [NS68]). In fact, if $m$ is a 
large integer which is a multiple of the denominators involved in
$G$ and $\Delta$, we can consider sections
$$\sigma\in H^0(X,m(K_X+G+\Delta)).$$
We view them rather as sections of $m(K_X+G)$ with poles along the 
support~$|\Delta|$ of our divisor. Then 
$(\sigma\wedge \overline\sigma)^{1/m}h_G$
is a volume form with
integrable poles along~$|\Delta|$ (this is the klt condition for~$\Delta$).
Therefore one can normalize $\sigma$ by requiring that
$$
\int_X (\sigma\wedge \overline\sigma)^{1/m}h_G=1.
$$
Each of these sections defines a singular hermitian metric on 
$K_X+L=K_X+G+\Delta$, and we can take the regularized upper envelope 
$$
\varphi_\can^\alg=\bigg(\sup_{m,\sigma}{1\over m}\log|\sigma|_{h_{K_X+L}^m}^2
\bigg)^*\leqno(5.18)
$$
of the weights associated with a smooth metric~$h_{K_X+L}$. It is clear that
$\varphi_\can^\alg\le \varphi_\can$ since the supremum is
taken on the smaller set of weights $\varphi={1\over m}
\smash{\log|\sigma|^2_{h_{K_X+L}^m}}$, and the equalities
$$
e^{\varphi-\gamma}dV_\omega=|\sigma|^{2/m}_{h_{K_X+L}^m}
e^{-\gamma}dV_\omega=(\sigma\wedge\ol\sigma)^{1/m}e^{-\gamma}h_L
=(\sigma\wedge\ol\sigma)^{1/m}h_{L,\gamma}
=(\sigma\wedge\ol\sigma)^{1/m}h_G
$$
imply $\int_X e^{\varphi-\gamma}dV_\omega\le 1$. We~claim that the inequality
$\varphi_\can^\alg\le \varphi_\can$ is an equality. The proof is an 
immediate consequence of the following
statement based in turn on the Ohsawa-Takegoshi theorem and
the approximation technique of [Dem92].}
\endclaim

\vbox{%
\claim (5.19) Proposition|With $L=G+\Delta$, $\omega$, 
$\alpha=\smash{\Theta_{K_X+L,h_{K_X+L}}}$, $\gamma$ 
as above and $K_X+L$ assumed to be~big, fix a singular hermitian metric 
$e^{-\varphi}h_{K_X+L}$ of curvature $\alpha+dd^c\varphi\ge 0$, such that
$\int_X e^{\varphi-\gamma}dV_\omega\le 1$. Then
$\varphi$ is equal to a regularized limit
\vskip-5pt
$$
\varphi=\bigg(\limsup_{m\to+\infty}{1\over m}\log|\sigma_m|_{h_{K_X+L}^m}^2
\bigg)^*
$$ 
for a suitable sequence $\sigma_m\in H^0(X,m(K_X+G+\Delta))$ with
$\int_X(\sigma_m\wedge\overline\sigma_m)^{1/m}h_G\le 1$.
\endclaim}

\proof. By our assumption, there exists a quasi-psh function $\psi_0$ with
analytic singularity set $Z_0$ such that
$$
\alpha+dd^c\psi_0\ge \varepsilon_0\omega>0
$$
and we can assume $\int_Ce^{\psi_0-\gamma}dV_\omega<1$ (the strict inequality
will be useful later). For
$m\ge p\ge 1$, this defines a singular metric 
$\exp(-(m-p)\varphi-p\psi_0)h_{K_X+L}^m$
on $m(K_X+L)$ with curvature${}\ge p\varepsilon_0\omega$, and therefore
a singular metric 
$$
h_{L'}=\exp(-(m-p)\varphi-p\psi_0)h_{K_X+L}^mh_{K_X}^{-1}
$$
on $L'=(m-1)K_X+mL$, whose curvature
$\Theta_{L',h_{L'}}\ge (p\varepsilon_0-C_0)\omega$ is arbitrary large
if $p$ is large enough. Let us fix a finite covering of $X$
by coordinate balls. Pick a point $x_0$ and one of the coordinate
balls $B$ containing~$x_0$. By the Ohsawa-Takegoshi extension theorem
applied on the ball~$B$,
we can find a section $\sigma_B$ of $K_X+L'=m(K_X+L)$
which has norm $1$ at $x_0$ with respect to the metric
$h_{K_X+L'}$ and $\int_B|\sigma_B|_{h_{K_X+L'}}^2dV_\omega\le C_1$ for some
uniform constant $C_1$ depending on the finite covering, but independent of 
$m$, $p$, $x_0$ . Now, we use a cut-off function $\theta(x)$ with
$\theta(x)=1$ near $x_0$ to truncate
$\sigma_B$ and solve a $\dbar$-equation for $(n,1)$-forms with values
in $L$ to get a global section $\sigma$ on $X$ with 
$|\sigma(x_0)|_{h_{K_X+L'}}=1$. For this we need to multiply our metric
by a truncated factor $\exp(-2n\theta(x)\log|x-x_0|)$ so as to get
solutions of $\dbar$ vanishing at~$x_0$. However, this perturbs the
curvature by bounded terms and we can absorb them again by taking
$p$ larger. In~this way we obtain
$$
\int_X|\sigma|^2_{h_{K_X+L'}}dV_\omega=
\int_X |\sigma|^2_{h_{K_X+L}^m}e^{-(m-p)\varphi-p\psi_0}dV_\omega\le C_2.
\leqno(5.20)
$$
Taking $p>1$, the H\"older inequality for congugate exponents
$m$, ${m\over m-1}$ implies
$$
\eqalign{
\int_X(\sigma\wedge\overline\sigma)^{1\over m}h_G
&=\int_X|\sigma|^{2/m}_{h_{K_X+L}^m}e^{-\gamma}dV_\omega\cr
&=\int_X\Big(|\sigma|^2_{h_{K_X+L}^m}e^{-(m-p)\varphi-p\psi_0}\Big)^{1\over m}
\Big(e^{(1-{p \over m})\varphi+{p\over m}\psi_0-\gamma}\Big)dV_\omega\cr
&\le C_2^{1\over m}\bigg(\int_X
\Big(e^{(1-{p\over m})\varphi+{p\over m}\psi_0-\gamma}
\Big)^{m\over m-1}dV_\omega\bigg)^{m-1\over m}\cr
&\le C_2^{1\over m}\bigg(\int_X
 \big(e^{\varphi-\gamma}\big)^{m-p\over m-1}
\Big(e^{{p\over p-1}(\psi_0-\gamma)}\Big)^{p-1\over m-1}dV_\omega
\bigg)^{m-1\over m}\cr
&\le C_2^{1\over m}
\bigg(\int_X e^{{p\over p-1}(\psi_0-\gamma)}dV_\omega
\bigg)^{p-1\over m}\cr}
$$
using the hypothesis $\int_X e^{\varphi-\gamma}dV_\omega\le 1$ and another
application of H\"older's inequality. Since 
klt is an open condition and $\lim_{p\to+\infty}
\smash{\int_Xe^{{p\over p-1}(\psi_0-\gamma)}dV_\omega}
=\int_Xe^{\psi_0-\gamma}dV_\omega<1$, we can take
$p$ large enough to ensure that 
$$
\int_X e^{{p\over p-1}(\psi_0-\gamma)}dV_\omega\le C_3<1.
$$
Therefore, we see that 
$$
\int_X(\sigma\wedge\overline\sigma)^{1\over m}h_G\le 
C_2^{1\over m}C_3^{p-1\over m}\le 1
$$
for $p$ large enough. On the other hand
$$
|\sigma(x_0)|^2_{h_{K_X+L'}}
=|\sigma(x_0)|^2_{h_{K_X+L}^m}e^{-(m-p)\varphi(x_0)-p\psi_0(x_0)}=1,
$$
thus
$$
{1\over m}\log |\sigma(x_0)|^2_{h_{K_X+L}^m}=
\Big(1-{p\over m}\Big)\varphi(x_0)+{p\over m}\psi_0(x_0)
\leqno(5.21)
$$
and, as a consequence
$$
{1\over m}\log |\sigma(x_0)|^2_{h_{K_X+L}^m}\longrightarrow\varphi(x_0)
$$
whenever $m\to+\infty$, ${p\over m}\to 0$, as long as
$\psi_0(x_0)>-\infty$. In the above argument, we can in fact interpolate
in finitely many points $x_1,\,x_2,\,\ldots\,,x_q$ provided that
$p\ge C_4q$. Therefore if we take a suitable dense subset
$\{x_q\}$ and a ``diagonal'' sequence associated with sections
$\sigma_m\in H^0(X,m(K_X+L))$ with 
$m\gg p=p_m\gg q=q_m\to+\infty$, we infer that
$$
\bigg(\limsup_{m\to+\infty}{1\over m}\log |\sigma_m(x)|^2_{h_{K_X+L}^m}\bigg)^*\ge
\limsup_{x_q\to x}\;\varphi(x_q)=\varphi(x)\leqno(5.22)
$$
(the latter equality occurring if $\{x_q\}$ is suitably chosen with
respect to $\varphi$). In the other direction, (5.20)
implies a mean value estimate
$$
{1\over \pi^nr^{2n}/n!}\int_{B(x,r)}
|\sigma(z)|^2_{h_{K_X+L}^m}\;dz\le {C_5\over r^{2n}}
\sup_{B(x,r)}e^{(m-p)\varphi+p\psi_0}
$$
on every coordinate ball $B(x,r)\subset X$. The function 
$|\sigma_m|^2_{h_{K_X+L}^m}$ is 
plurisubharmonic after we correct the
non necessarily positively curved smooth metric $h_{K_X+L}$ by a 
factor of the form $\exp(C_6|z-x|^2)$, hence the mean value 
inequality shows that
$$
{1\over m}\log |\sigma_m(x)|^2_{h_{K_X+L}^m}\le 
{1\over m}\log{C_5\over r^{2n}}+C_6 r^2
+\sup_{B(x,r)}\Big(1-{p_m\over m}\Big)\varphi+{p_m\over m}\psi_0.
$$
By taking in particular $r=1/m$ and letting $m\to+\infty$, $p_m/m\to 0$, 
we see that the opposite of inequality (5.22) also holds.\qed

\claim (5.23) Remark|{\rm We can rephrase our results in slightly different
terms. In fact, let us put 
$$
\varphi^\alg_m=\sup_{\sigma}{1\over m}\log|\sigma|_{h_{K_X+L}^m}^2,
\qquad \sigma\in H^0(X,m(K_X+G+\Delta)),
$$
with normalized sections $\sigma$
such that \hbox{$\int_X (\sigma\wedge \overline\sigma)^{1/m}h_G=1$}.
Then $\varphi^\alg_m$ is quasi-psh (the supremum is taken over a
compact set in a finite dimensional vector space) and by passing 
to the regularized supremum over all $\sigma$ and all $\varphi$ 
in (5.21) we get
$$
\varphi_\can\ge\varphi^\alg_m\ge \Big(1-{p\over m}\Big)
\varphi_\can(x)+{p\over m}\psi_0(x).
$$
As $\varphi_\can$ is bounded from above, we find in particular
$$
0\le \varphi_\can-\varphi^\alg_m\le {C\over m}(|\psi_0(x)|+1).
$$
This implies that $(\varphi^\alg_m)$ converges uniformly to
$\varphi_\can$ on every compact subset of $X\subset Z_0$, and in this way
we infer again (in a purely qualitative manner) that $\varphi_\can$ 
is continuous on $X\ssm Z_0$. Moreover, we also see that in (5.18) 
the upper  semicontinuous regularization is not needed on
$X\ssm Z_0\,$; in case $K_X+L$ is ample, it is not needed at all
and we have uniform convergence of $(\varphi^\alg_m)$ towards
$\varphi_\can$ on the whole of~$X$. Obtaining such a uniform convergence
when $K_X+L$ is just big looks like a more delicate question, 
related e.g.\ to abundance
of $K_X+L$ on those subvarieties $Y$ where the restriction
$(K_X+L)_{|Y}$ would be e.g.\ nef but not big.}
\endclaim

\claim (5.24) Generalization|{\rm In the general case
where $L$ is a $\bR$-line bundle and $K_X+L$ is merely pseudo-effective, 
a similar algebraic approximation can be obtained. We take instead sections
$$
\sigma\in H^0(X,mK_X+\lfloor mG\rfloor+\lfloor m\Delta\rfloor +p_mA)
$$
where $(A,h_A)$ is a positive line bundle, 
$\Theta_{A,h_A}\ge \varepsilon_0\omega$,
and replace the definition of $\varphi_\can^\alg$ by
$$
\leqalignno{
&\varphi_\can^\alg=
\bigg(\limsup_{m\to+\infty}\;
\sup_\sigma\;{1\over m}\log|\sigma|_{h_{mK_X+\lfloor mG\rfloor+p_mA}}^2
\bigg)^*,&(5.25)\cr
&\int_X(\sigma\wedge\overline\sigma)^{2\over m}
h_{\lfloor mG\rfloor +p_mA}^{1\over m}\le 1,&(5.26)\cr}
$$
where $m\gg p_m\gg 1$ and $h_{\lfloor mG\rfloor}^{1/m}$ is chosen
to converge uniformly to $h_G$.

 We then find again $\varphi_\can=\varphi_\can^\alg$, with
an almost identical proof -- though we no longer have a $\sup$ in the
envelope, but just a $\limsup$. The analogue of Proposition (5.19) also
holds true in this context, with an appropriate sequence of sections
$\sigma_m\in H^0(X,mK_X+\lfloor mG\rfloor+\lfloor m\Delta\rfloor +p_mA)$.}
\endclaim

\claim (5.27) Remark|{\rm The envelopes considered in section 1 are
envelopes constrained by an $L^\infty$ condition, while the
present ones are constrained by an $L^1$ condition. It is 
possible to interpolate and to consider envelopes constrained by
an $L^p$ condition. More precisely, assuming that
${1\over p}K_X+L$ is pseudo-effective, we look at metrics
$e^{-\varphi}h_{{1\over p}K_X+L}$ and normalize them with the $L^p$ condition
$$
\int_X e^{p\varphi-\gamma}dV_\omega\le 1.
$$
This is actually an $L^1$ condition for the induced metric on $pL$, and 
therefore we can just apply the above after replacing $L$ by $pL$. 
If we assume moreover that
$L$ is pseudo-effective, it is clear that the $L^p$ condition converges
to the $L^\infty$ condition $\varphi\le 0$, if we normalize $\gamma$ by
requiring that $\int_X e^{-\gamma}dV_\omega=1$.}
\endclaim

\claim (5.28) Remark|{\rm It would be nice to have a better understanding
of the supercanonical metrics. In case $X$ is a curve, this should be easier.
In fact $X$ then has a hermitian metric $\omega$ with constant curvature,
which we normalize by requiring that $\int_X\omega=1$, and we can also
suppose $\int_Xe^{-\gamma}\omega=1$. The class
$\lambda=c_1(K_X+L)\ge 0$ is a number and we take $\alpha=\lambda\omega$. 
Our envelope is $\varphi_\can=\sup\varphi$ where 
$\lambda\omega+dd^c\varphi\ge 0$ and
$\int_Xe^{\varphi-\gamma}\omega\le 1$. If $\lambda=0$ then $\varphi$
must be constant and clearly $\varphi_\can=0$. \hbox{Otherwise}, if $G(z,a)$ 
denotes the Green function such that $\int_X G(z,a)\,\omega(z)=0$ and
$dd^cG(z,a)=\delta_a-\omega(z)$, we find
$$
\varphi_\can(z)\ge \sup_{a\in X}\bigg(\lambda
G(z,a)-\log\int_{z\in X}e^{\lambda G(z,a)-\gamma(z)}\omega(z)\bigg)
$$
by taking already the envelope over $\varphi(z)=\lambda G(z,a)-\hbox{Const}$. 
It is natural
to ask whether this is always an equality, i.e.\ whether the extremal
functions are always given by one of the Green functions, especially
when $\gamma=0$.
}
\endclaim

\section{References}
\medskip

{\eightpoint

\bibitem[BT76]&Bedford, E.\ {\rm and} Taylor, B.A.:& The Dirichlet 
  problem for the complex Monge-Amp\`ere equation;& Invent.\ Math.\
  {\bf 37} (1976) 1--44&

\bibitem[BT82]&Bedford, E.\ {\rm and} Taylor, B.A.:& A new capacity for 
  plurisubharmonic functions;& Acta Math.\ {\bf 149} (1982) 1--41&

\bibitem[Ber07]&Berman, R.:& Bergman kernels and equilibrium measures for
  line bundles over projective manifolds;& arXiv:0710.4375,
  to appear in the American Journal of Mathematics&

\bibitem[BB08]&Berman, R., Boucksom, S.:& Growth of balls of holomorphic
  sections and energy at equilibrium;& arXiv:0803.1950&

\bibitem[BBGZ09]&Berman, R., Boucksom, S., Guedj, V., Zeriahi, A.:& 
  A variational approach to solving complex Monge-Amp\`ere equations;&
  in preparation&

\bibitem[Blo09]&B{\l}ocki, Z:& On geodesics in the space of K\"ahler 
  metrics;& Preprint 2009 available at\\ {\tt
  http://gamma.im.uj.edu.pl/\~{}blocki/publ/}&

\bibitem[Bou02]&Boucksom, S.:& On the volume of a line bundle;& Internat.\
  J.\ Math {\bf 13} (2002), 1043--1063&

\bibitem[BEGZ08]&Boucksom, S., Eyssidieux, P., Guedj, V., Zeriahi, A.:&
  Monge-Amp\`ere equations in big cohomology classes;& arXiv:0812.3674&

\bibitem [CNS86]&Caffarelli, L., Nirenberg, L., Spruck, J.:& The Dirichlet 
  problem for the degenerate Monge-Amp\`ere equation;& Rev.\ Mat.\ 
  Iberoamericana {\bf 2} (1986), 19--27&

\bibitem [Che00]&Chen, X.:&The space of K\"ahler metrics;& J.\ Differential
  Geom.\ {\bf 56} (2000), 189--234&

\bibitem[CT08]&Chen, X. X., Tian, G.:& Geometry of K\"ahler metrics and
  foliations by holomorphic discs;& Publ.\ Math.\ Inst.\ Hautes \'Etudes 
  Sci.\ {\bf 107} (2008), 1--107&

\bibitem[Dem89]&Demailly, J.-P.:& Potential Theory in 
  Several Complex Variables;&
  Manuscript available at 
  www-fourier.ujf-grenoble.fr/$\tilde{~}$demailly/books.html, [D3]&

\bibitem[Dem91]&Demailly, J.-P.:& Complex analytic and algebraic geometry;&
  manuscript Institut Fourier, first edition 1991, available online at\\
  {\tt http://www-fourier.ujf-grenoble.fr/\~{}demailly/books.html} &

\bibitem[Dem92]&Demailly, J.-P.:& Regularization of closed positive currents 
  and Intersection Theory;& J.\ Alg.\ Geom.\ {\bf 1} (1992), 361--409&

\bibitem[Dem93]&Demailly, J.-P.:& Monge-Amp\`ere operators, Lelong numbers 
  and intersection theory;& Complex Analysis and Geometry, Univ.\ Series in 
  Math., edited by V.~Ancona and A.~Silva, Plenum Press, New-York, 1993&

\bibitem[Dem94]&Demailly, J.-P.:& Regularization of closed positive 
  currents of type $(1,1)$ by the flow of a Chern connection;& 
  Actes du Colloque en l'honneur de P.~Dolbeault
  (Juin 1992), \'edit\'e par H.~Skoda et J.-M.~Tr\'epreau,
  Aspects of Mathematics, Vol.~E~26, Vieweg, 1994, 105--126&

\bibitem [Don99]&Donaldson, S.K.:& Symmetric spaces, K\"ahler geometry and 
  Hamiltonian dynamics;& Northern California Symplectic Geometry Seminar, 
  13--33, Amer.\ Math.\ Soc.\ Transl.\ Ser.\ 2, 196, 
  Amer.\ Math.\ Soc., Providence, RI, 1999&

\bibitem [Don02]&Donaldson, S.K.:& Holomorphic discs and the complex
  Monge-Amp\`ere equation;& J.\ Symplectic Geom.\ {\bf 1} (2002), 171--196&

\bibitem[Kis78]&Kiselman, C.~O.:& The partial Legendre transformation
  for plurisubharmonic functions;& Inventiones Math.\ {\bf 49} (1978)
  137--148&

\bibitem[Kis94]&Kiselman, C.~O.:& Attenuating the singularities of
  plurisubharmonic functions;& Ann.\ Polonici Mathematici {\bf 60} (1994)
  173--197&

\bibitem[NS68]&Narasimhan, M.S., Simha, R.R.:& Manifolds with ample canonical
  class;& Inventiones Math.\ {\bf 5} (1968) 120--128&

\bibitem[OhT87]&Ohsawa, T., Takegoshi, K.:& On the extension
  of $L^2$ holomorphic functions;& Math.\ Zeitschrift {\bf 195} (1987)
  197--204&

\bibitem[PS08]&Phong, D.H., Sturm. J:& The Dirichlet problem for 
  degenerate complex Monge-Amp\`ere equations;& arXiv:0904.1898&

\bibitem[Rud66]&Rudin, W.:& Real and complex analysis;& 
  McGraw-Hill, 1966, third edition 1987&

\bibitem[Ts07a]&Tsuji, H.:& Canonical singular hermitian metrics on
  relative canonical bundles;& arXiv:0704.0566&

\bibitem[Ts07b]&Tsuji, H.:& Canonical volume forms on compact K\"ahler
  manifolds;& arXiv:0707.0111&

}
\vskip6pt
\noindent
(version of May 5, 2009, printed on \today)

\end